\documentclass[12pt]{amsart}
\usepackage{amssymb,amscd}
\usepackage[all]{xy}
\title[Arakelov (in)equalities]{Arakelov (in)equalities}
\author[E. Viehweg]{Eckart Viehweg}
\address{Universit\"at Duisburg-Essen, Mathematik, 45117 Essen, Germany}
\email{viehweg@uni-due.de}
\thanks{This work has been supported by the DFG-Leibniz program and by the SFB/TR 45
``Periods, moduli spaces and arithmetic of algebraic varieties''.}
 \fussy
\begin{document}
\theoremstyle{plain}
\newtheorem{thm}{Theorem}[section]
\newtheorem{theorem}[thm]{Theorem}
\newtheorem{addendum}[thm]{Addendum}
\newtheorem{variant}[thm]{Variant}
\newtheorem{lemma}[thm]{Lemma}
\newtheorem{corollary}[thm]{Corollary}
\newtheorem{proposition}[thm]{Proposition}
\newtheorem{some_results}[thm]{Some Results}
\theoremstyle{definition}
\newtheorem{remark}[thm]{Remark}
\newtheorem{remarks}[thm]{Remarks}
\newtheorem{notations}[thm]{Notations}
\newtheorem{definition}[thm]{Definition}
\newtheorem{claim}[thm]{Claim}
\newtheorem{assumption}[thm]{Assumption}
\newtheorem{assumptions}[thm]{Assumptions}
\newtheorem{property}[thm]{Property}
\newtheorem{properties}[thm]{Properties}
\newtheorem{example}[thm]{Example}
\newtheorem{examples}[thm]{Examples}
\newtheorem{conjecture}[thm]{Conjecture}
\newtheorem{questions}[thm]{Questions}
\newtheorem{question}[thm]{Question}
\newtheorem{construction}[thm]{Construction}
\numberwithin{equation}{section}
\newcommand{\sA}{{\mathcal A}}
\newcommand{\sB}{{\mathcal B}}
\newcommand{\sC}{{\mathcal C}}
\newcommand{\sD}{{\mathcal D}}
\newcommand{\sE}{{\mathcal E}}
\newcommand{\sF}{{\mathcal F}}
\newcommand{\sG}{{\mathcal G}}
\newcommand{\sH}{{\mathcal H}}
\newcommand{\sI}{{\mathcal I}}
\newcommand{\sJ}{{\mathcal J}}
\newcommand{\sK}{{\mathcal K}}
\newcommand{\sL}{{\mathcal L}}
\newcommand{\sM}{{\mathcal M}}
\newcommand{\sN}{{\mathcal N}}
\newcommand{\sO}{{\mathcal O}}
\newcommand{\sP}{{\mathcal P}}
\newcommand{\sQ}{{\mathcal Q}}
\newcommand{\sR}{{\mathcal R}}
\newcommand{\sS}{{\mathcal S}}
\newcommand{\sT}{{\mathcal T}}
\newcommand{\sU}{{\mathcal U}}
\newcommand{\sV}{{\mathcal V}}
\newcommand{\sW}{{\mathcal W}}
\newcommand{\sX}{{\mathcal X}}
\newcommand{\sY}{{\mathcal Y}}
\newcommand{\sZ}{{\mathcal Z}}
\newcommand{\A}{{\mathbb A}}
\newcommand{\B}{{\mathbb B}}
\newcommand{\C}{{\mathbb C}}
\newcommand{\D}{{\mathbb D}}
\newcommand{\E}{{\mathbb E}}
\newcommand{\F}{{\mathbb F}}
\newcommand{\G}{{\mathbb G}}
\newcommand{\BH}{{\mathbb H}}
\newcommand{\I}{{\mathbb I}}
\newcommand{\J}{{\mathbb J}}
\newcommand{\BL}{{\mathbb L}}
\newcommand{\M}{{\mathbb M}}
\newcommand{\N}{{\mathbb N}}
\newcommand{\BP}{{\mathbb P}}
\newcommand{\Q}{{\mathbb Q}}
\newcommand{\R}{{\mathbb R}}
\newcommand{\BS}{{\mathbb S}}
\newcommand{\T}{{\mathbb T}}
\newcommand{\U}{{\mathbb U}}
\newcommand{\V}{{\mathbb V}}
\newcommand{\W}{{\mathbb W}}
\newcommand{\X}{{\mathbb X}}
\newcommand{\Y}{{\mathbb Y}}
\newcommand{\Z}{{\mathbb Z}}
\newcommand{\rk}{{\rm rk}}
\newcommand{\ch}{{\rm c}}
\newcommand{\Sp}{{\rm Sp}}
\newcommand{\Sl}{{\rm Sl}}
\newcommand{\bsA}{{\overline{\sA}}}
\newcommand{\BEnd}{{\mathbb E}{\rm nd}}
\newcommand{\Mon}{{\rm Mon}}
\newcommand{\MT}{{\rm MT}}
\newcommand{\Hg}{{\rm Hg}}
\newcommand{\End}{{\rm End}}
\newcommand{\Gal}{{\rm Gal}}
\newcommand{\doubledot}{{\mbox{\Tiny\textbullet,\textbullet}}}
\newcommand{\ds}{\displaystyle}
\newcommand{\SU}{{\rm SU}}
\maketitle
\section*{Introduction}
The proof of the Shafarevich Conjecture for curves of genus $g\geq 2$ over complex function fields $K=\C(Y)$, given by Arakelov in \cite{Ar71}, consists of two parts, the verification of ``boundedness'' and of ``rigidity''.
In order to obtain the boundedness, Arakelov first constructs a height function for $K$-valued points of the moduli stack
$\overline{\sM}_g$ of stable curves of genus $g$. In down to earth terms, he chooses a natural ample sheaf $\lambda$ on the coarse moduli scheme $\overline{M}_g$. Then, extending the morphism ${\rm Spec}(K) \to \overline{M}_g$
to $Y \to \overline{M}_g$ he chooses as height $\deg(\varphi^*\lambda)$.
Secondly, still assuming that $\varphi$ is induced by a genuine family $f:X\to Y$ of stable curves, he gives an upper bound for this height in terms of the curve $Y$ and the discriminant $S=Y\setminus Y_0$ for $Y_0= \varphi^{-1}(M_g)$. Finally the rigidity, saying that $X_0=f^{-1}(Y_0) \to Y_0$ does not extend to a family ${\mathfrak f}:{\mathfrak X}_0\to Y_0\times T$ in a non-trivial way,
easily follows from the deformation theory for families of curves. 

The boundedness part of Arakelov's proof was extended by Faltings \cite{Fa83} to families of abelian varieties, using
Deligne's description of abelian varieties via Hodge structures of weight one. He chooses on a suitable toroidal compactification
$\overline{\sA}_g$ of the coarse moduli scheme of polarized abelian varieties and $\lambda\in {\rm Pic}(\overline{\sA}_g)\otimes \Q$ to be the determinant of the direct image of relative one forms, hence the determinant of the Hodge bundle of bidegree $(1,0)$ in the corresponding variation of Hodge structures. Then $\lambda$ is semiample and ample with respect to the open set
$\sA_g$ (as defined in Definition \ref{1.2}), which is sufficient to define a height function.
He proves an upper bound for the height, hence the finiteness of deformation types, and gives a criterion for infinitesimal rigidity. A family of $8$-dimensional abelian varieties gives an example that contrary to the case of curves the rigidity fails in general.  

Deligne \cite{De87} takes up Faltings approach. He obtains more precise inequalities and his arguments extend to
$\C$-variations of Hodge structures of weight one. Peters proved similar inequalities for variation of Hodge structures of higher weight. Unfortunately his results (improved by Deligne in an unpublished letter) were 
only available years later (see \cite{Pe00}), shortly after the subject was taken up by Jost and Zuo in \cite{JZ02}.  

Since then the results for families of curves or abelian varieties over curves have been extended in several ways.  
Firstly the definition and the bounds for height functions have been extended to moduli schemes of canonically polarized manifolds
or of polarized minimal models (see \cite{BV00}, \cite{vz01}, \cite{vz04a}, \cite{V05} and \cite{kl06}, for example). 
We sketch some of the results in Section \ref{s.1}. However we will not say anything about rigidity and strong boundedness properties, discussed in \cite{vz02} and \cite{kl06}. 

Secondly generalizations of the Arakelov inequalities are known for variations of Hodge structures of higher weight over curves, and for weight one over a higher dimensional bases. In both cases the inequalities are optimal, i.e. there are families where one gets equality. As we recall in Section \ref{s.1} such an equality should be rare for families of varieties of positive Kodaira dimension. Except for abelian varieties and for K3-surfaces the geometric interpretation of such an equality is still not understood.  

Finally the Arakelov inequalities have a topological counterpart, the Milnor-Wood inequalities for the Toledo invariant, for certain local systems on projective curves and on higher dimensional projective manifolds (see \cite{bgg06}, \cite{km08a}, and \cite{km08b}, for example) . Again the equality has consequences for the structure of the local system (or its Higgs bundle). We will state this (in)equalities in very special cases in Section \ref{s.5} and in Section \ref{s.8} and compare it with the Arakelov inequality.

The main theme of this survey is the interplay between stability of Higgs bundles and the stability of the Hodge bundles for variations of Hodge structures of weight $k$ (see Section \ref{s.2} for the basic definitions). As we try to explain in Section \ref{s.3} for all $k$ in the curve case, and in Section \ref{s.6} for $k=1$ over certain higher dimensional varieties, the Arakelov inequalities are translations of slope conditions for polystable Higgs bundles, whereas the Arakelov equalities encode stability conditions for the Hodge bundles. In Sections \ref{s.4} and \ref{s.7} we indicate some geometric consequences of Arakelov equalities for $k=1$ or for families of abelian varieties.\\[.2cm]
{\bf Acknowledgments.}
This survey is based on a series of articles coauthored by Kang Zuo, by Martin M\"oller or by both of them. Compared with those articles there are only minor improvements in some arguments and no new results.

Martin M\"oller pointed out some ambiguities in the first version of this article, and 
the idea for the simplified proof of Claim~\ref{6.7}, needed for Theorem~\ref{6.4}, is taken from his letter
explaining the ``$r=2$''-case. I am gratefull to Vincent Koziarz and Julien Maubon for their explanations
concerning the ``Milnor-Wood'' inequality over a higher dimensional base.

\section{Families of manifolds of positive Kodaira dimension}\label{s.1}
Let $f:X\to Y$ be a semistable family of $n$-folds over a complex projective curve $Y$, smooth over $Y_0=Y\setminus S$ and with $X$ projective. We call $f$ {\em semistable} if $X$ is non-singular and if all fibres $f^{-1}(y)$ of $f$ are reduced normal crossing divisors. We write $X_0=f^{-1}(Y_0)$ and $f_0=f|_{X_0}$. 
\begin{theorem}[\cite{vz01}, \cite{vz06}, and \cite{mvz06}]\label{1.1}
Assume that $f:X\to Y$ is semistable. Then for all $\nu\geq 1$ with $f_*\omega_{X/Y}^\nu\neq 0$
\begin{equation}\label{eq1.1}
\frac{\deg (f_* \omega^{\nu}_{X/Y})}{\rk(f_*
\omega^{\nu}_{X/Y})}\leq \frac{n\cdot \nu}{2} \cdot \deg(\Omega^1_Y(\log S)).
\end{equation}
\end{theorem}
The morphism $f$ is called {\em isotrivial} if there is a finite covering $Y'\to Y$ and a birational $Y'$ morphism 
$$X\times_YY' \dashrightarrow F\times Y'.$$ 
For projective manifolds $F$ with $\omega_F$ semiample and polarized by an invertible sheaf with Hilbert polynomial $h$, there exists a coarse quasiprojective moduli scheme $M_h$. Hence if $\omega_{X_0/Y_0}$ is $f_0$-semiample
$f_0$ induces a morphism $\varphi_0:Y_0\to M_h$.

If $\omega_{X_0/Y_0}$ is $f_0$-ample, or if $\omega_{X_0/Y_0}^\nu$ is for some $\nu>0$ the pullback of an invertible sheaf on $Y_0$, then the birational non-isotriviality of $f$ is equivalent to the quasi-finiteness of $\varphi_0$. In this situation   
the left hand side of ~\eqref{eq1.1} can be seen as a height function on the moduli scheme.
In fact, choosing $\nu>1$ with $h(\nu)\neq 0$ in the first case, and or $\nu\geq 1$ with $\omega_F^\nu=\sO_F$ in the second one, by \cite{V05} there exists a projective compactification $\overline{M}_h$ of the moduli scheme
$M_h$ and some 
$$
\lambda\in {\rm Pic}(\overline{M}_h)\otimes \Q
$$ 
with:
\begin{itemize}
\item $\lambda$ is nef and ample with respect to $M_h$.
\item Let $\varphi:Y\to \overline{M}_h$ be the morphism induced by $f$. Then
$\det(f_* \omega^{\nu}_{X/Y})=\varphi^*\lambda$.
\end{itemize}
For moduli of abelian varieties one can choose the Baily-Borel compactification and there $\lambda$ is ample. 
By \cite{mu77} on a suitable toroidal compactification of $\sA_g$ the sheaf $\lambda$ is still semi-ample,
but for other moduli functors we only get weaker properties, as defined below.
\begin{definition}\label{1.2} Let $Z$ be a projective variety and let $Z_0\subset Z$ be open and dense.
\begin{enumerate}
\item[i.] A locally free sheaf $\sF$ on $Z$ is {\em numerically effective (nef)} if for all morphisms $\rho:C\to Z$, with $C$ an irreducible curve, and for all invertible quotients $\sN$ of $\rho^*\sF$ one has $\deg(\sN) \geq 0$.
\item[ii.] An invertible sheaf $\sL$ on $Z$ is {\em ample with respect to $Z_0$} if for some $\nu\geq 1$ the sections in $H^0(Z,\sL^\nu)$ generate the sheaf $\sL^\nu$ over $Z_0$ and if the induced morphism $Z_0 \to \BP(H^0(Z,\sL^\nu))$ is an embedding.  
\end{enumerate}
\end{definition}
For non-constant morphisms $\rho:C\to Z$ from irreducible projective curves one finds in Definition~\ref{1.2}, ii) that $\deg(\rho^*(\sL))>0$, provided $\rho(C)\cap Z_0\neq \emptyset$. Moreover, fixing an upper bound $c$ for this degree, there are only finitely many deformation types of curves with $\deg(\rho^*(\sL))<c$. 

Applying this to birationally non-isotrivial families $f:X\to Y$ whose general fibre $F$ is either canonically polarized or a minimal model of Kodaira dimension zero, one finds the left hand side of~\eqref{eq1.1} to be positive, hence $\Omega^1_Y(\log S)=\omega_Y(S)$ must be ample. The finiteness of the number of deformation types is more difficult and it has been
worked out in \cite{kl06} just for families of canonically polarized manifolds. 
Roughly speaking, one has to show that morphisms from a curve to the moduli stack are parameterized by a scheme. This being done, one finds that for a given Hilbert polynomial $h$ and for a given constant $c$ there are only finitely many deformation types of families $f:X\to Y$ of canonically polarized manifolds with $\deg(\Omega^1_Y(\log S))\leq c$.

For smooth projective families $f_0:X_0\to Y_0$ over a higher dimensional quasi-projective manifold $Y_0$ with $\omega_{X_0/Y_0}$ semiample, some generalizations of the inequality~\eqref{eq1.1} have been studied in \cite{vz02} (see also \cite{vz04a}). There we assumed that $S=Y\setminus Y_0$ is a normal crossing divisor and that the induced map $\varphi_0:Y_0\to M_h$ is generically finite. Then for some $\mu \gg 0$ there exists a non-trivial ample subsheaf of $S^\mu(\Omega_Y^1(\log S))$. However neither $\mu$ nor the degree of the ample subsheaf have been calculated and the statement is less precise than the inequality~\eqref{eq1.1}.

In this survey we are mainly interested in a geometric interpretation of equality in~\eqref{eq1.1}, in particular for $\nu=1$. As explained in \cite{vz06} and \cite{mvz06} such equalities should not occur for families with $p_g(F)>1$ for the general fibre $F$. Even the Arakelov inequalities for non-unitary subvariations of Hodge structures, discussed in Section~\ref{s.3} should be strict for {\em most} families with $F$ of general type. As recalled in Example \ref{4.6}, for curves ``most'' implies that the genus $g$ of $F$ has to be $3$ and that the ``counter-example'' in genus 3 is essentially unique. So what Arakelov equalities are  concerned it seems reasonable to concentrate on families of minimal models of Kodaira dimension zero.
 
\section{Stability}\label{s.2}

\begin{definition}\label{2.1} Let $Y$ be a projective manifold, let $S\in Y$ be a normal crossing divisor and let $\sF$ be a torsion-free coherent sheaf on $Y$.
\begin{enumerate}
 \item[i.] The degree and slope of $\sF$ are defined as
$$
\deg(\sF) = \ch_1(\sF).\ch_1({\omega_Y(S)})^{\dim(Y)-1} \mbox{ \ \ and \ \ }  \mu(\sF)=\mu_{\omega_Y(S)}(\sF)=\frac{\deg(\sF)}{\rk(\sF)}. 
$$
\item[ii.] The sheaf $\sF$ is $\mu$-stable if for all subsheaves $\sG\subset \sF$
with $\rk(\sG)<\rk(\sF)$ one has $\mu(\sG)<\mu(\sF)$.\vspace{.1cm}
\item[iii.] The sheaf $\sF$ is $\mu$-semistable if for all non-trivial subsheaves $\sG\subset \sF$
one has $\mu(\sG)\leq\mu(\sF)$.\vspace{.1cm}
\item[iv.] $\sF$ is $\mu$-polystable if it is the direct sum of $\mu$-stable sheaves of the same slope. 
\end{enumerate}
\end{definition}
This definition is only reasonable if $\dim(Y)=1$ or if $\omega_Y(S)$ is nef and big.

Recall that a {\em logarithmic Higgs bundle} is a locally free sheaf $E$ on $Y$ together with an $\sO_Y$ linear morphism $\theta:E\to E\otimes \Omega^1_Y(\log S)$ with $\theta\wedge\theta=0$. 
The definition of stability (poly- and semistability) for locally free sheaves extends to Higgs bundles, by requiring that
$$
\mu(F)=\frac{\deg(F)}{\rk(F)} < \mu(E)=\frac{\deg(E)}{\rk(E)} 
$$
(or $\mu(F)\leq \mu(E)$) for all subsheaves $F$ with $\theta(F)\subset F\otimes \Omega^1_Y(\log S)$.

If $\dim(Y) > 1$, for the Simpson correspondence in \cite{Si92} and for the polystability of Higgs bundles, one 
takes the slopes with respect to a polarization of $Y$, i.e replacing
$\omega_Y(S)$ in Definition~\ref{2.1}, i) by an ample invertible sheaf. 
However, as we will recall in Proposition~\ref{6.4}, the Simpson correspondence remains true
for the slopes $\mu(\sF)$ in~\ref{2.1}, i), provided $\omega_Y(S)$ is nef and big. 

Our main example of a Higgs bundle will be the one attached to a 
polarized $\C$ variation of Hodge structures $\V$ on $Y_0$ of weight $k$, as defined in \cite{De87}, and with unipotent local monodromy operators. The $\sF$-filtration of $F_0=\V\otimes_\C \sO_{Y_0}$ extends to a locally splitting filtration of the Deligne extension $F$ of $F_0$ to $Y$, denoted here by
$$
\sF^{k+1} \subset \sF^{k} \subset \cdots \subset \sF^0.
$$
We will usually assume that $\sF^{k+1}=0$ and $\sF^0=F$, hence that all non-zero parts of the Hodge decomposition of a fibre $\V_y$ of $\V$ are in bidegrees $(k-m,m)$ for $m=0,\ldots, k$.
The Griffiths transversality condition for the Gau{\ss}-Manin connection
$\nabla$ says that
$$
\nabla(\sF^p)\subset \sF^{p-1}\otimes \Omega^1_Y(\log S).
$$
Then $\nabla$ induces a $\sO_Y$ linear map
$$
\theta_{p,k-p}: E^{p,k-p}=\sF^p/\sF^{p+1} \longrightarrow E^{p-1,k-p+1}=\sF^{p-1}/\sF^{p}\otimes
\Omega^1_Y(\log S).
$$
We will call
$$
\big(E=\bigoplus_p E^{p,k-p}, \theta=\bigoplus \theta_{p,k-p}\big)
$$
the {\em (logarithmic) Higgs bundle of $\V$}, whereas the sheaves $E^{p,q}$ are called the
{\em Hodge bundles} of bidegree $(p,q)$.
\begin{definition}\label{2.2} For the Higgs bundle $(E,\theta)$ introduced above we define: 
\begin{enumerate}
\item[i.] The {\em support} ${\rm supp}(E,\theta)$ is the set of all $m$ with $E^{k-m,m}\neq 0$.
\item[ii.] $(E,\theta)$ has a {\em connected support}, if there exists some $m_0 \leq m_1 \in \Z$ with
\begin{gather*}
{\rm supp}(E,\theta)=\{m;\ m_0 \leq m \leq m_1\} \mbox{ \ \ and if \ \ }\\
\theta_{k-m,m}\neq 0 \mbox{ \ \ for \ \ }
m_0\leq m \leq m_1-1.
\end{gather*}
\item[iii.] $(E,\theta)$ (or $\V$) satisfies the {\em Arakelov condition}
if $(E,\theta)$ has a connected support and if for all $m$ with $m, \ m+1 \in {\rm supp}(E,\theta)$ the sheaves $E^{k-m,m}$
and $E^{k-m-1,m+1}$ are $\mu$-semistable and 
$$
\mu(E^{k-m,m})=\mu(E^{k-m-1,m+1}) + \mu(\Omega^1_Y(\log S)).
$$
\end{enumerate}
\end{definition}
\section{Variations of Hodge structures over curves}\label{s.3}
Let us return to a projective curve $Y$, so $S=Y\setminus Y_0$ is a finite set of points. The starting point of our considerations is the Simpson correspondence:
\begin{theorem}[\cite{Si90}]\label{3.1}
There exists a natural equivalence between the category of direct
sums of stable filtered regular Higgs bundles of degree zero, and
of direct sums of stable filtered local systems of degree zero.
\end{theorem}
We will not recall the definition of a ``filtered regular'' Higgs bundle
\cite[page 717]{Si90}, and just remark that for a Higgs bundle corresponding
to a local system $\V$ with unipotent monodromy around the points in
$S$ the filtration is trivial, and automatically $\deg(\V)=0$.

By \cite{De71} the local systems underlying a $\Z$-variation of Hodge structures are semisimple, and by
\cite{De87} the same holds with $\Z$ replaced by $\C$. So one obtains:
\begin{corollary}\label{3.2}
The logarithmic Higgs bundle of a polarized $\C$-variation of Hodge structures with unipotent monodromy in $s\in S$ is polystable
of degree $0$.
\end{corollary}
In \cite{vz03} and \cite{vz06} we discussed several versions of Arakelov inequalities. Here we will only need the one for $E^{k,0}$, and we sketch a simplified version of the proof:
\begin{lemma}\label{3.3} Let $\V$ be an irreducible complex polarized variation of Hodge structures over $Y$
of weight $k$ and with unipotent local monodromies in $s\in S$.
Write $(E,\theta)$ for the logarithmic Higgs bundle of $\V$ and assume that $E^{p,k-p}=0$ for
$p<0$ and for $p>k$. Then one has:\vspace{.1cm}
\begin{enumerate}
\item[a.] \ \hspace*{\fill} \ \ \ \ \ \ \ \ \ $\displaystyle 
\mu(E^{k,0}) \leq \frac{k}{2} 
\cdot \deg(\Omega^1_Y(\log S)).$ \hspace*{\fill} \ \vspace{.1cm}
\item[b.] \ \hspace*{\fill} $ \ \  0 \leq \mu(E^{k,0}) $ \hspace*{\fill} \vspace{.1cm} \\
and the equality implies that $\V$ is unitary or equivalently that $\theta=0$.\vspace{.1cm}
\item[c.] The equality \ \hspace*{2.2cm} $\displaystyle 
\mu(E^{k,0}) = \frac{k}{2} 
\cdot \deg(\Omega^1_Y(\log S)).$ \hspace*{\fill} \vspace{.1cm} \\
implies that the sheaves $E^{k-m,m}$ are  stable and that 
$$
\theta_{k-m,m}:E^{k-m,m}\longrightarrow E^{k-m-1,m+1}\otimes \Omega^1_Y(\log S)
$$ 
is an isomorphism for $m=0, \ldots, k-1$.   
\end{enumerate}
\end{lemma}
\begin{proof} Let $G^{k,0}$ be a subsheaf of $E^{k,0}$, and 
let $G^{k-m,m}$ be the $(k-m,m)$ component of the Higgs subbundle $G=\langle G^{k,0} \rangle$, generated by $G^{k,0}$. By definition one has a surjection
$$
G^{k-m+1,m-1} \longrightarrow G^{k-m,m} \otimes \Omega^1_Y(\log S).
$$ 
Its kernel $\sK_{m-1}$, together with the $0$-map is a Higgs subbundle of $(E,\theta)$, hence of non-positive degree.
Remark that $\rk(G^{k,0}) \geq \rk(G^{k-1,1}) \geq \cdots \geq \rk(G^{k-m,m})$. So one finds
\begin{multline}\label{eq3.1}
\deg(G^{k-m+1,m-1}) \leq \deg(G^{k-m,m}) + \rk(G^{k-m,m})\cdot \deg(\Omega^1_Y(\log S))\leq \\
\deg(G^{k-m,m}) + \rk(G^{k-1,1})\cdot \deg(\Omega^1_Y(\log S)).
\end{multline}
Iterating this inequality gives for $m\geq 1$
\begin{multline}\label{eq3.2}
\deg(G^{k,0}) \leq \deg(G^{k,0}) - \deg(\sK_0) = \\
\deg(G^{k-1,1}) + \rk(G^{k-1,1})\cdot \deg(\Omega^1_Y(\log S)) \leq \\
\deg(G^{k-m,m}) + m \cdot \rk(G^{k-1,1})\cdot \deg(\Omega^1_Y(\log S))
\end{multline}
and adding up 
\begin{multline*}
(k+1) \deg(G^{k,0}) \leq (k+1) \deg(G^{k,0}) - k \cdot \deg(\sK_0) \leq\\ \sum_{m=0}^{k} \deg(G^{k-m,m})
+  \sum_{m=1}^{k} m \cdot \rk(G^{k-1,1})\cdot \deg(\Omega^1_Y(\log S))=\\
\deg(G) + \frac{k\cdot(k+1)}{2}\cdot \rk(G^{k-1,1})\cdot \deg(\Omega^1_Y(\log S)).
\end{multline*}
Since $G$ is a Higgs subbundle, $\deg(G) \leq 0$, and 
\begin{equation}\label{eq3.3} 
\mu (G^{k,0}) \leq \frac{\deg(G^{k,0})}{\rk(G^{k-1,1})} \leq \frac{k}{2} \cdot \deg(\Omega^1_Y(\log S)).
\end{equation}
Taking $G^{k,0}=E^{k,0}$ one obtains the inequality in a). 

If this is an equality, as assumed in c), then the right hand side of \eqref{eq3.3} is an equality.
Firstly, since the difference of the two sides is larger than a positive multiple of $\deg(G)=0$, the latter is zero and the irreducibility of $\V$ implies that $G=E$. Secondly the two inequalities  
in \eqref{eq3.2} have to be equalities. The one on the right hand side gives $\rk(E^{k-m,m})=\rk(E^{k-1,1})$ for $m=2,\ldots, k$. The one on the left implies that $\deg(\sK_0)=0$ and the irreducibility of $\V$ shows that this is only possible for $\sK_0=0$ hence if $\rk(E^{k,0})=\rk(E^{k-1,1})$. All together one finds that the surjections 
$$
E^{k,0} \longrightarrow E^{k-m,m} \otimes \Omega^1_Y(\log S)^m
$$ 
are isomorphisms, for $1\leq m \leq k$. On the other hand the equality in c) and the inequality
\eqref{eq3.3} imply that for all subsheaves $G^{k,0}$ 
$$
\mu (G^{k,0}) \leq \frac{k}{2} \cdot \deg(\Omega^1_Y(\log S))= \mu (E^{k,0}),
$$
If this is an equality, then $\deg(G)=0$ and $(G,\theta|_G)\subset (E,\theta)$ splits. The irreducibility implies again that
$(G,\theta|_G)=(E,\theta)$, hence $E^{k,0}$ as well as all the $E^{k-m,m}$ are  stable.

The sheaf $E^{k,0}$ with the $0$-Higgs field is a Higgs quotientbundle of $(E,\theta)$, hence of non-negative degree.
If $\deg(E^{k,0})=0$, then the surjection of Higgs bundles $(E,\theta)\to (E^{k,0},0)$ splits. The irreducibility
of $\V$ together with Theorem~\ref{3.1} implies that both Higgs bundles are the same, hence that $\theta=0$ and $\V$ unitary. So b) follows from a).
\end{proof}
\begin{corollary}\label{3.4}
In Lemma \ref{3.3} one has the inequality
\begin{equation}\label{eq3.4} 
\mu(E^{k,0}) - \mu(E^{0,k}) \leq k \cdot \deg(\Omega^1_Y(\log S)).
\end{equation}
The equality in Lemma \ref{3.3}, c) is equivalent to the equality
\begin{equation}\label{eq3.5}
\mu(E^{k,0}) - \mu(E^{0,k}) = k \cdot \deg(\Omega^1_Y(\log S)).
\end{equation}
In particular \eqref{eq3.5} implies that the sheaves $E^{k-m,m}$ are  stable and that 
$$
\theta_{k-m,m}:E^{k-m,m}\longrightarrow E^{k-m-1,m+1}\otimes \Omega^1_Y(\log S_)
$$ 
is an isomorphism for $m=0, \ldots, k-1$.
\end{corollary}
\begin{proof}
For \eqref{eq3.4} one applies part a) of Lemma \ref{3.3} to $(E,\theta)$ and to the dual Higgs bundle $(E^\vee,\theta^\vee)$. 
The equality \eqref{eq3.5} implies that both, $(E,\theta)$ and $(E^\vee,\theta^\vee)$ satisfy the Arakelov equality
c) in Lemma \ref{3.3}. 

Finally assume that the equation c) in Lemma \ref{3.3} holds for $(E,\theta)$. Then 
\begin{gather*}
E^{k,0}\cong E^{0,k}\otimes \Omega^1_Y(\log S)^k \mbox{ \ \ and}\\
\mu({E^\vee}^{k,0})=-\mu(E^{0,k})=  k\cdot \deg(\Omega^1_Y(\log S)) -\mu(E^{k,0}) = \frac{k}{2} \cdot \deg(\Omega^1_Y(\log S)).
\end{gather*}
Adding this equality to the one in c) one gets \eqref{eq3.5}.
\end{proof}

The inequality in part a) of Lemma \ref{3.3} is not optimal. One can use the degrees of the kernels $\sK_m$ to get correction terms. We will only work this out for $m=0$.  What equalities are concerned, one does not seem to get anything new. 
\begin{variant}\label{3.5}
In Lemma~\ref{3.6} one has the inequalities
\begin{equation}\label{eq3.6}
\frac{\deg(E^{k,0})}{\rk(\theta_{k,0})} \leq \frac{k}{2} \cdot \deg(\Omega^1_Y(\log S)).
\end{equation}
The equality in Lemma \ref{3.3}, c) is equivalent to the equality
\begin{equation}\label{eq3.7}
\frac{\deg(E^{k,0})}{\rk(\theta_{k,0})} = \frac{k}{2} \cdot \deg(\Omega^1_Y(\log S)).
\end{equation}
\end{variant}
\begin{proof}
The inequality is a repetition of the left hand side of~\eqref{eq3.3} for $G^{k,0} = E^{k,0}$.
 
If $\theta_{k,0}$ is an isomorphisms, hence if $\rk(E^{k,0})=\rk(\theta_{k,0})$, the two
equalities~\eqref{eq3.7} and c) in Lemma \ref{3.3} are the same. As stated in Lemma \ref{3.3}, the equality c)
implies that $\theta_{k,0}$ is an isomorphisms, hence \eqref{eq3.7}.

In the proof of Lemma \ref{3.3} we have seen that the equality of the right hand side of \eqref{eq3.3} implies that $G=E$, hence that the morphisms 
$$
\theta_{k-m,m}:E^{k-m,m} \longrightarrow E^{k-m-1,m+1}\otimes \Omega^1_Y(\log S)
$$ 
are surjective for $m=0,\ldots,m-1$. Using the left hand side of \eqref{eq3.2}, one finds that $\sK_0=0$ hence that $\theta_{k,0}$ is an isomorphisms. So \eqref{eq3.7} implies the equality c). 
\end{proof}

Replacing $Y_0$ by an \'etale covering, if necessary, one may assume that $\# S$ is even, hence that there exists 
a logarithmic theta characteristic $\sL$. By definition $\sL^2\cong \Omega^1_Y(\log S)$ and one has an isomorphism
$$
\tau: \sL \longrightarrow \sL\otimes \Omega^1_Y(\log S).
$$
Since $(\sL\oplus \sL^{-1},\tau)$ is an indecomposable Higgs bundle of degree zero, Theorem~\ref{3.1}
tell us that it comes from a local system $\BL$, which is easily seen to be a variation of Hodge structures of weight $1$. We will say that $\BL$ is induced by a logarithmic theta characteristic. Remark that $\BL$ is unique
up to the tensor product with local systems, corresponding to two division points in ${\rm Pic}^0(Y)$.
By \cite[Proposition 3.4]{vz03} one has:
\begin{addendum}\label{3.6} Assume in Lemma~\ref{3.3} that $\# S$ is even and that $\BL$ is induced by a theta characteristic. 
\begin{enumerate}
\item[d.] Then the equality $\displaystyle 
\mu(E^{k,0}) = \frac{k}{2} 
\cdot \deg(\Omega^1_Y(\log S))$ implies that there exists an irreducible unitary local system $\T_0$ on $Y_0$ with
$\V \cong \T_0 \otimes S^k(\BL).$
\end{enumerate}     
\end{addendum}
\begin{remark}\label{3.7}
In Addendum \ref{3.6} the local monodromies of $\T_0$ are unipotent and unitary, hence finite. So there exists a finite covering $\tau:Y'\to Y$, \'etale over $Y_0$ such that $\tau^*\T_0$ extends to a unitary local system $\T'$ on $Y'$.
\end{remark}
The property d) in Addendum \ref{3.6} is equivalent to the condition c) in Lemma \ref{3.3}. In particular it implies that each $E^{k-m,m}$ is the tensor product of an invertible sheaf with the polystable sheaf $\T_0\otimes_\C\sO_Y$. The Arakelov equality implies that the Higgs fields are direct sums of morphisms between semistable sheaves of the same slope. Then the irreducibility of $\V$ can be used to show that $\T_0\otimes_\C\sO_Y$ and hence the $E^{k-m,m}$ are stable.
\begin{remark} \label{3.8} 
Let us collect what we learned in the proof of Lemma~\ref{3.3}. 
\begin{itemize}
\item Simpson's polystability of the Higgs bundles $(E,\theta)$
implies the Arakelov inequality a) in Lemma \ref{3.3} or inequality \eqref{eq3.4}. 
\item The equality in part c) of Lemma~\ref{3.3} implies that the Hodge bundles $E^{k-m,m}$ are semistable and that the Higgs field is a morphism of sheaves of the same slope.
\item If one assumes in addition that $\V$ is irreducible, then the $E^{k-m,m}$ are stable sheaves.
\end{itemize}
As we will see in Section \ref{s.6} the first two statements extend to families over a higher dimensional base (satisfying the positivity condition ($\star$) in \ref{6.2}), but we doubt that the third one remains true without some additional numerically conditions.
\end{remark}
Assume that $\W$ is the variation of Hodge structures given by a smooth family $f_0:X_0\to Y_0$ of polarized manifolds with semistable reduction at infinity, hence $\W=R^kf_{0*}\C_{X_0}$. Let $\W=\V_1\oplus \cdots \oplus \V_\ell$ be the decomposition of $\W$ as direct sum of irreducible local subsystems, hence of $\C$ irreducible variations of Hodge structures of weight $k$. Replacing $\V_\iota$ by a suitable Tate twist $\V_\iota(\nu_\iota)$, and perhaps by its dual, one obtains a variation of Hodge structures of weight $k_\iota=k-2\cdot \nu_\iota$, whose Hodge bundles
are concentrated in bidegrees $(k_\iota-m,m)$ for $m=0,\ldots,k_\iota$ and non-zero in bidegree
$(k_\iota,0)$. Applying Lemma \ref{3.3} to $\V_\iota(\nu_\iota)$ one gets Arakelov inequalities for
all the $\V_\iota$. If all those are equalities, each of the $\V_\iota$ will satisfy the Arakelov condition
in Definition \ref{2.2}, iii, and for some unitary  bundle $\T_\iota$ one finds
$\V_\iota=\T_\iota \otimes S^{k-2\cdot\nu_\iota}(\BL)(-\nu_\iota)$. We say that the Higgs field of $\W$ is {\em strictly maximal} in this case (see \cite{vz03} for a motivation and for a slightly different presentation of those results).

Let us list two results known for families of Calabi-Yau manifolds, satisfying the Arakelov equality. 
\begin{assumptions}\label{3.9}
Consider smooth morphisms $f_0:X_0\to Y_0$ over a non-singular curve $Y_0$, whose fibres are $k$-dimensional Calabi-Yau manifolds. Assume that $f_0$ extends to a semistable family $f:X\to Y$ on the compactification $Y$ of $Y_0$.
Let $\V$ be the irreducible direct factor of $R^kf_{0*}\C_{X_0}$ with Higgs bundle $(E,\theta)$, such that $E^{k,0}\neq 0$. 
\end{assumptions}
\begin{theorem}[\cite{Bo97}, \cite{Vo93}, and \cite{STZ03}, see also \cite{vz03}]\label{3.10}\ \\
For all $k\geq 1$ there exist families $f_0:X_0\to Y_0$ satisfying the Assumptions \ref{3.9}, such that the Arakelov equality \eqref{eq3.5} holds for $\V$. For families of $K3$-surfaces, i.e. for $k=2$, there exist examples with $Y_0=Y$ projective.
\end{theorem}\label{3.11}
For $k=1$ those families are the universal families over elliptic modular curves, hence $Y_0$ is affine in this case. A similar result holds whenever the dimension of the fibres is odd. 
\begin{theorem}[\cite{vz03}] Under the assumptions made in \ref{3.9} assume that $k$ is odd and that
 $\V$ satisfies the Arakelov equality. Then $S=Y\setminus Y_0 \neq 0$, i.e. $Y_0$ is affine.
\end{theorem}
It does not seem to be known whether for even $k\geq 4$ there are  families of Calabi-Yau manifolds over a compact curve with $\V$ satisfying the Arakelov equality. 

The geometric implications of the Arakelov equality for $\V$ in \ref{3.9} or of the
strict maximality of the Higgs field, are not really understood. The structure Theorem \ref{3.6} can be used to obtain some properties of the Mumford Tate group, but we have no idea about the structure of the family or about  the map to the moduli scheme $M_h$. The situation is better for families of abelian varieties.
So starting from the next section we will concentrate on polarized variations of Hodge structures of weight one. 

\section{Arakelov equality and geodecity of curves in $\sA_g$}\label{s.4}
\begin{assumptions}\label{4.1}
Keeping the assumptions from the last section, we restrict ourselves to variations of Hodge structures of weight one, coming from families $f_0:X_0\to Y_0$ of abelian varieties. Replacing $Y_0$ by an \'etale covering allows to assume that
$f_0:X_0\to Y_0$ is induced by a morphism $\varphi_0:Y_0\to \sA_g$ where $\sA_g$ is some fine moduli scheme of polarized abelian varieties with a suitable level structure, and that the local monodromy in $s\in S$ of
$\W_\Q=R^1f_{0*}\Q_{X_0}$ is unipotent. Let us fix a toroidal compactification $\bsA_g$, as considered by Mumford in \cite{mu77}. In particular $\bsA_g$ is non-singular, the boundary divisor $S_{\bsA_g}$ has non-singular components, and normal crossings, $\Omega_{\bsA_g}^1(\log S_{\bsA_g})$ is nef and $\omega_{\bsA_g}(S_{\bsA_g})$
is ample with respect to $\sA_g$.
\end{assumptions}
In \cite{mv08} we give a differential geometric characterization of morphisms $\varphi_0:Y_0\to \sA_g$ 
for which the induced $\C$-variation of Hodge structures $\W$ contains a non-unitary
$\C$-subvariation $\V$ with Higgs bundle $(E,\theta)$, satisfying the Arakelov equality
\begin{equation}\label{eq4.1}
\mu (E^{1,0}) = \frac{1}{2} \cdot \deg(\Omega^1_Y(\log S)).
\end{equation}
To this aim we need:
\begin{definition}\label{4.2} Let $M$ be a complex domain and $W$ be a subdomain. $W$
is a {\em totally geodesic submanifold for the Kobayashi metric} if 
the restriction of the Kobayashi metric on $M$ to $W$ coincides with the
Kobayashi metric on $W$. If $W = \Delta$ we call $\Delta$
{\em a (complex) Kobayashi geodesic}.
\par 
A map $\varphi_0: Y_0 \to A_g$ is {\em a Kobayashi geodesic}, if
its universal covering map 
$$
\tilde{\varphi}_0: \widetilde{Y_0} \cong \Delta \longrightarrow \BH_g
$$ is a Kobayashi geodesic. In particular here a Kobayashi geodesic will always be
one-dimensional.
\end{definition}
\begin{theorem}\label{4.3} Under the assumptions made in~\ref{4.1} the following conditions are equivalent:
\begin{enumerate}
\item[a.] $\varphi_0: Y_0 \to \sA_g$ is Kobayashi geodesic.
\item[b.] The natural map $\varphi^* \Omega^1_{\bsA_g}(\log S_{\bsA_g}) \to \Omega^1_Y(\log S)$
splits.
\item[c.] $\W$ contains a non-unitary irreducible subvariation of Hodge structures $\V$ which satisfies the Arakelov equality \eqref{eq4.1}.
\end{enumerate}
\end{theorem}
The numerical condition in Theorem~\ref{4.3} indicates that  Kobayashi geodesic in $\sA_g$ are
``algebraic objects''. In fact, as shown in  \cite{mv08} one obtains:
\begin{corollary}\label{4.4}
Let $\varphi_0:Y_0\to \sA_g$ be an affine Kobayashi geodesic, such that the
induced variation of Hodge structures $\W_\Q$ is $\Q$-irreducible.
Then $\varphi_0:Y_0 \to \sA_g$ can be defined over a number field.   
\end{corollary}
Geodesics for the Kobayashi metric have been considered in \cite{Moe06} under the additional assumption that $f_0:X_0 \to Y_0$ is a family of Jacobians of a smooth family of curves. In this case $\varphi_0(Y_0)$ is a geodesic for the Kobayashi metric if and only if the image of $Y_0$ in the moduli scheme $M_g$ of curves of genus $g$ with the right level structure is a geodesic for the Teichm\"uller metric, hence if and only if $Y_0$ is a Teichm\"uller curve. In particular $Y_0$ will be affine and the 
irreducible subvariation $\V$ in Theorem \ref{4.3} will be of rank two. By Addendum \ref{3.6} it is given by a logarithmic theta characteristic on $Y$. Using the theory of Teichm\"uller curves (see \cite{McM03}), one can deduce that there is at most one irreducible direct factor $\V$ which satisfies the Arakelov equality.

The Theorem~\ref{4.3} should be compared with the results of \cite{vz04}. 
Starting from Lemma~\ref{3.3} and the addendum~\ref{3.6} it is shown that under the assumptions~\ref{4.1}
$Y_0$ (or to be more precise, an \'etale finite cover of $Y_0$) is a rigid Shimura curve with universal family $f_0:X_0\to Y_0$ if the Arakelov equality holds for all irreducible $\C$-subvariations of Hodge structures of $R^1f_{0*}\C_{X_0}$. Recall that ``rigid'' means that there are no non-trivial extensions of $f_0$ to a smooth family ${\mathfrak f}: {\mathfrak X}_0 \to T\times Y_0$ with $\dim{T}>0$.
If one allows unitary direct factors, and requires the Arakelov equality just for all non-unitary
subvariations $\V$, then $Y_0 \subset \sA_g$ is a deformation of a Shimura curve or, using the notation from
\cite{Mu69}, the family $f_0:X_0\to Y_0$ is a Kuga fibre space.  

In \cite{Moe05} it is shown (see also \cite[Section 1]{mvz07}), that for all Kuga fibre spaces
and all non-unitary irreducible $\V\subset R^1f_{0*}\C_{X_0}$ the Arakelov equality holds.
In \cite{mvz07} this was translated to geodecity for the Hodge (or Bergman-Siegel) metric, and we can restate
the main result of \cite{vz04} in the following form:
\begin{theorem}\label{4.5}
Keeping the notations and assumptions introduced in \ref{4.1}, the following conditions are equivalent:
\begin{enumerate}
\item[a.] $\varphi_0: Y_0 \to \sA_g$ is a geodesic for the Hodge metric on $\sA_g$.
\item[b.] The natural map $\varphi^* \Omega^1_{\bsA_g}(\log S_{\bsA_g}) \to \Omega^1_Y(\log S)$
splits orthogonal for the Hodge metric.
\item[c.] All non-unitary irreducible $\C$-subvariations of Hodge structures $\V \subset \W$ satisfy the Arakelov equality.
\item[d.] $f_0:X_0\to Y_0$ is a Kuga fibre space over the curve $Y_0$.
\end{enumerate} 
\end{theorem}
\begin{example}\label{4.6}
One can ask, whether a geodesic for the  Hodge metric on $\sA_g$ can lie completely in $M_g$, or in different terms,
whether there exists a family of smooth curves over $Y_0$, such that the induced family of Jacobians is a Kuga fibre space.
This is of course true for modular families of elliptic curves.   

In \cite{vz06} it was shown that for a Hodge geodesic in $M_g$ the rang of the maximal non-unitary part of the corresponding variation of Hodge structures has to be $2$. So by \cite{Moe06} $Y_0$ is a Teichm\"uller curve and $Y_0$ is affine.   
By \cite{Moe05} the only ``Shimura-Teichm\"uller curve'', i.e. the only Hodge geodesics in $M_g$, exists for  $g=3$. Up to \'etale coverings, there is only one example.
\end{example}
\section{Milnor-Wood inequalities}
\label{s.5}
Before we discuss families of abelian varieties over a higher dimensional base, let us mention a
numerical condition, which applies to a different class of Higgs bundles over curves $Y$, the 
Milnor-Wood inequality for the Toledo invariant. We refere to \cite{bgg06} for an introduction and for a guide to the literature.

Let $\T$ be a local system, induced by a representation of $\pi_1(Y,*)$ in a connected non-compact semi-simple real Lie group $G$. Since the representations of the fundamental group of $Y$ are not semi-simple we can not apply Simpson's correspondence stated in Theorem \ref{3.1}. As in \cite[Section 2]{bgg06} one has to add on the representation side the condition ``reductive'' and on the Higgs bundle side 
the condition ``polystable''. 

As explained in \cite[Section 3.2]{bgg06}, for $G={\rm SU}(p,q)$ (or for $G=\Sp(2n,\R)$) the corresponding Higgs bundle $(F,\psi)$ is given by two locally free sheaves $\sV$ and $\sW$ on $Y$ of rank $p$ and $q$, respectively, and the Higgs field $\psi$ is the direct sum of two morphisms 
$$
\beta: \sW \longrightarrow \sV \otimes \Omega_Y^1 \mbox{ \ \ and \ \ }
\gamma: \sV \longrightarrow \sW \otimes \Omega_Y^1
$$
(For $G=\Sp(2n,\R)$ one has $p=q=n$ and $\sW=\sV^\vee$). 
The Toledo invariant of $\T$ or of $(F,\tau)$ is 
$$
\tau(\T)=\tau((F,\psi))=\deg(\sV)=-\deg(\sW), 
$$
and the classical Milnor-Wood inequality says that
$$
|\tau(\T)| \leq {\rm Min}\{p,q\}\cdot (g-1).
$$
In fact, on page 194 of \cite{bgg06} one finds a more precise inequality, and again equality has strong implications on the structure of the Higgs field:

\begin{proposition}\label{5.1}
Let $\T$ be a local system on $Y_0$ induced by a representation of $\pi_1(Y,*)$ in ${\rm SU}(p,q)$ (or in $\Sp(2n,\R)$). Assume that the Higgs bundle 
$$
(F=\sV\oplus \sW,\psi=\gamma +\beta)
$$ 
is polystable, where 
$$
\gamma: \sV  \longrightarrow \sW\otimes \Omega_Y^1(\log S)
\mbox{ \ \ and \ \ } \beta: \sW \longrightarrow \sV\otimes \Omega_Y^1(\log S).
$$
Then 
\begin{equation}\label{eq5.1}
-\rk(\beta)\cdot \deg(\Omega^1_Y)
 \leq - 2\cdot \deg(\sW)=2\cdot \deg(\sV) \leq \rk(\gamma) \cdot \deg(\Omega^1_Y).
\end{equation}
The inequality $2\cdot \deg(\sV) \leq \rk(\gamma) \cdot \deg(\Omega^1_Y)$ is strict, except if
$\gamma$ is an isomorphism.  
\end{proposition}
\begin{proof}
It is sufficient to prove the inequality on the right hand side. The other one follows by interchanging the role of $p$ and $q$, hence of $\sV$ and $\sW$. Since this inequality is compatible with exact sequences of Higgs bundles, the  Jordan-H\"older filtration for Higgs bundles allows to assume that $(F=\sV\oplus \sW,\psi=\gamma+\beta)$ is stable.

The subbundle $\sG=\sV\oplus \gamma(\sV)\otimes (\Omega^1_Y)^{-1}$ of $(F,\theta)$ and the kernel $\sK$ of
$\sV \to \gamma(\sV)$ are compatible with the Higgs field, hence
\begin{multline*}
2\cdot\deg(\sV) - \rk(\gamma(\sV))\cdot \deg(\Omega^1_Y) \leq\\
\deg(\sV) + \deg(\gamma(\sV)) - \rk(\gamma(\sV))\cdot \deg(\Omega^1_Y) \leq 0.
\end{multline*}
If equality holds, the stability implies that $\sK$ is zero and that $\sG=F$. In particular $\gamma$ is an isomorphism.
\end{proof}
\begin{example}\label{5.2}
Let $(E,\theta)$ be the Higgs bundle of a polarized variation of Hodge structures of weight one.
Then one could choose $\sV=E^{1,0}$ and $\sW=E^{0,1}$. Since $\beta=0$ and $\theta=\gamma$
the inequality~\eqref{eq5.1} says that
$$
0 \leq 2 \cdot \deg(\sV) \leq \rk(\gamma) \cdot \deg(\Omega^1_Y),
$$ 
hence it coincides with the inequality~\eqref{eq3.6} in Variant~\ref{3.5}.
\end{example}

\begin{example}[Kang Zuo]\label{5.3}
Let $(E,\theta)$ be the Higgs field of a variation of Hodge structures of weight $k$, with $k$ odd.
Choose
$$
\sV=\bigoplus_{m=0}^{\frac{k-1}{2}} E^{k-2m,2m} \mbox{ \ \ and \ \ }
\sW=\bigoplus_{m=0}^{\frac{k-1}{2}} E^{k-2m-1,2m+1},
$$
and for $\gamma$ and $\beta$ the restriction of the Higgs field. The Milnor-Wood inequality
says that 
\begin{equation}\label{eq5.2}
-\sum_{m=0}^{\frac{k-1}{2}} \deg(E^{k-2m-1,2m+1}) = \sum_{m=0}^{\frac{k-1}{2}} \deg(E^{k-2m,2m}) \leq \frac{1}{2}\cdot \deg(\Omega^1_Y).
\end{equation}
The Arakelov equality in Lemma \ref{3.3}, c) implies that \eqref{eq5.2} is an equality. On the other hand,
having equality in \eqref{eq5.2} just implies that the morphisms 
$$
\theta_{k-m,m}:E^{k-m,m} \longrightarrow E^{k-m-1,m+1}\otimes \Omega^1_Y 
$$
are isomorphisms for $m$ even, but it says nothing about the other components of the Higgs field.
So the two equalities are not equivalent. To get an explicit example, consider
$Y=\BP^1$ and $S=\{0,1,\infty\}$. Choose 
$$
E^{3,0}=\sO_{\BP^1}(1), \ \ E^{2,1}=\sO_{\BP^1} \ \ E^{1,2}=\sO_{\BP^1} \ \ E^{0,3}=\sO_{\BP^1}(-1).
$$
So $\theta_{3,0}$ and $\theta_{1,2}$ are isomorphisms, whereas $\theta_{2,1}:\sO_{\BP^1}\to  \Omega_{\BP^1}^1(\log (0+1+\infty))$
is injective and has a zero in some point, say $2$. The degree of $E$ is zero, and obviously the Higgs bundle is stable. Hence by
Theorem \ref{3.1} $(E,\theta)$ is the Higgs bundle of a local system, and by construction the local system underlies a polarized variation of Hodge structures. 

In addition this example shows that \eqref{eq5.1} can be an equality, without $\sV$ or $\sW$ being stable.
\end{example}
\section{Arakelov inequalities for variations of Hodge structures of weight one over a higher dimensional base}
\label{s.6}
From now on $Y$ denotes a projective manifold and $S$ a normal crossing divisor in $Y$ with $Y_0=Y\setminus S$. We will need some positivity properties of the sheaf of differential forms on the compactification $Y$ of $Y_0$.
\begin{assumptions}\label{6.1}
We suppose that
\begin{enumerate}
\item[($\star$)] \ \hspace*{\fill} $\Omega_Y^1(\log S)$ is nef and $\omega_Y(S)$ is ample with respect to $Y_0$.\hspace*{\fill} \ 
\end{enumerate}
\end{assumptions}
As a motivation, assume for the moment that $Y_0 \subset \sA_g$ is a Shimura variety, or that there is a Kuga fibre space $f_0:X_0\to Y_0$. In both cases $Y_0$ is the quotient of a bounded symmetric domain and replacing $Y_0$ by an \'etale finite cover, we may choose for $Y$ a {\em Mumford compactification}, i.e. a toroidal compactification, as studied in \cite{mu77}. There it is shown that $\Omega^1_Y(\log S)$ is nef, and that $Y$ maps to the Baily-Borel compactification. Then the finiteness of $Y_0\to \sA_g$ implies that $\omega_Y(S)$ is ample with respect to $Y_0$. 

The main reason why we need an extra condition is Yau's Uniformization Theorem
(\cite{Ya93}, discussed in \cite[Theorem 1.4]{vz07}), saying that ($\star$) forces the sheaf $\Omega_Y^1(\log S)$ to be $\mu$-polystable. Here, as in Definition~\ref{2.1}, we will consider for coherent sheaves $\sF$ the slope $\mu(\sF)$ with respect to $\omega_Y(S)$. 

Usually to define stability and semistability on higher dimensional projective schemes, one considers
slopes with respect to polarizations. Replacing ``ample'' by ``nef and big'', hence considering {\em semi-polarizations} there might exist effective boundary divisors $D$ not recognized by the 
slope. So one has to identify $\mu$-equivalent subsheaves.
\begin{definition}\label{6.2} \ 
\begin{enumerate}
\item[1.] A subsheaf $\sG$ of $\sF$ is $\mu$-equivalent to $\sF$, if $\sF/\sG$ is a torsion sheaf and if $\ch_1(\sF)-\ch_1(\sG)$ is the class of an effective divisor $D$ with $\mu(\sO_Y(D))=0$.
\item[2.] $\sG\subset \sF$ is saturated, if $\sF/\sG$ is torsion free.
\item[3.] $\sF$ is weakly $\mu$-polystable, if it is $\mu$-equivalent to a $\mu$-polystable subsheaf.
\end{enumerate}
\end{definition}
The way it is stated, Theorem~\ref{3.1} only generalizes to a higher dimensional base $Y_0$ if $Y_0=Y$ is compact.
For variations of Hodge structures however the polystability of the induced Higgs bundle remains true and,  
as recalled in \cite[Proposition 2.4]{vz07}, there is no harm in working with the 
``semi-polarization'' $\omega_Y(S)$. 
\begin{proposition}\label{6.3} Let $E$ be the logarithmic Higgs bundle of a $\C$-variation of Hodge structures $\W$
on $Y_0$ with unipotent local monodromy around the components of $S$. If $G\subset E$ is a sub-Higgs sheaf then
for all $\dim(Y)-1\geq \nu \geq 0$ and for all ample invertible sheaves $\sH$ on $Y$ one has
\begin{equation}\label{eq6.1} 
\ch_1(G).\ch_1(\omega_Y(S))^{\dim(Y)-\nu-1}.\ch_1(\sH)^\nu \leq 0.
\end{equation}
Moreover, if $G\subset E$ is saturated the following conditions are equivalent:
\begin{enumerate}
\item For some $\nu \geq 0$ and for all ample invertible sheaves $\sH$ the equality holds in (\ref{eq6.1}).
\item For all $\nu $ and for all ample invertible sheaves $\sH$ the equality holds in (\ref{eq6.1}).
\item $G$ is induced by a local sub-system of $\W$.
\end{enumerate}
\end{proposition}
Most of the standard properties of stable and semistable sheaves carry over to the case of a semi-polarization, hence to the slope with respect to the nef and big sheaf $\omega_Y(S)$. In particular one can show the existence of maximal $\mu$-semistable destabilizing subsheaves, hence the existence of the Harder-Narasimhan filtration.
In addition the tensor product of $\mu$-semistable sheaves is again $\mu$-semistable. 

The starting point of \cite{vz07} was a generalization of the Arakelov inequality 
\eqref{eq3.5} for $k=1$ to a higher dimensional base. 
It is a direct consequence of the Simpson correspondence, as stated in Proposition~\ref{6.3}, using 
quite annoying calculations of slopes and degrees. Following a suggestion of Martin M\"oller we present below a simplified version of those calculations.
\begin{theorem}\label{6.4} Under the Assumption ($\star$) consider a polarized $\C$-variation of Hodge structures $\V$ on $Y_0$ with logarithmic Higgs bundle $(E,\theta)$. Assume that $\V$ is non-unitary, irreducible and that the local monodromy in $s\in S$ is unipotent. Then
\begin{equation}\label{eq6.2}
\mu(\V):= \mu(E^{1,0})-\mu(E^{0,1}) \leq \mu(\Omega^1_Y(\log S)). 
\end{equation} 
The equality $\mu(\V) = \mu(\Omega^1_Y(\log S))$ implies that
$E^{1,0}$ and $E^{0,1}$ are both $\mu$-semistable.
\end{theorem}
\begin{proof}
Consider the Harder-Narasimhan filtrations 
$$
0=G_0 \subsetneqq G_1 \subsetneqq \cdots \subsetneqq G_\ell=E^{1,0} \mbox{ \ \ and \ \ }
0=G'_0 \subsetneqq G'_1 \subsetneqq \cdots \subsetneqq G'_{\ell'}=E^{0,1}.
$$
Next choose two sequences of maximal length
\begin{gather}\notag
0=j_0 < j_1 < \cdots < j_r=\ell \mbox{ \ \ and \ \ } 0=j'_0 < j'_1 < \cdots < j'_r=\ell', \mbox{ \ \ with}\\ \label{eq6.3}
\theta(G_{j_{\iota}}) \subset G'_{j'_\iota}\otimes \Omega^1_Y(\log S), \hspace{1cm}
\theta(G_{j_{\iota-1}+1}) \not\subset G'_{j'_\iota-1}\otimes\Omega^1_Y(\log S)\\
\mbox{ and hence \ \ } \notag
\theta(G_{j_{\iota}}) \not\subset G'_{j'_\iota-1}\otimes\Omega^1_Y(\log S).
\end{gather}
Starting with $j_0=j'_0=0$ this can be done in the following way.
Assume one has defined $j_{\iota-1}$ and $j'_{\iota-1}$. Then
$j'_\iota$ is the minimal number with 
$$\theta(G_{j_{\iota-1}+1}) \subset G'_{j'_\iota}\otimes
\Omega^1_Y(\log S),$$
and $j_\iota$ is the maximum of all $j$ with
$$\theta(G_{j}) \subset G'_{j'_\iota}\otimes \Omega^1_Y(\log S).$$
Writing $E^{1,0}_\iota=G_{j_\iota}$ and $E^{0,1}_\iota=G'_{j'_\iota}$ we obtained two filtrations
$$
0=E^{1,0}_0 \subsetneqq E^{1,0}_1 \subsetneqq \cdots \subsetneqq E^{1,0}_r=E^{1,0} \mbox{ \ \ and \ \ }
0=E^{0,1}_0 \subsetneqq E^{0,1}_1 \subsetneqq \cdots \subsetneqq E^{0,1}_r=E^{0,1}.
$$
Let us define $F^{p,q}_\iota=E^{p,q}_\iota/E^{p,q}_{\iota-1}$. Remark that $F^{p,q}_\iota$ is not necessarily $\mu$-semistable. For $(p,q)=(1,0)$ for example, the Harder Narasimhan filtration is given by
$$
0=G_{j_{\iota-1}}/G_{j_{\iota-1}} \subsetneqq G_{j_{\iota-1}+1}/G_{j_{\iota-1}}  \subsetneqq \cdots \subsetneqq 
G_{j_\iota}/G_{j_{\iota-1}}=F^{1,0}_\iota.
$$
So 
\begin{equation}\label{eq6.4}
 \mu(G_{j_{\iota-1}+1}/G_{j_{\iota-1}}) \geq \mu(F^{1,0}_\iota) \geq \mu(G_{j_\iota}/G_{j_\iota-1})
\end{equation}
and, replacing $G$ by $G'$ and $j$ by $j'$,
\begin{equation}\label{eq6.5}
 \mu(G'_{j'_{\iota-1}+1}/G'_{j'_{\iota-1}}) \geq \mu(F^{0,1}_\iota) \geq \mu(G'_{j'_\iota}/G'_{j'_\iota-1}).
\end{equation}
\begin{claim}\label{6.5} \ 
\begin{enumerate}
\item[A.] \hspace*{0.5cm} $(\ch_1(E^{1,0}_\iota) + \ch_1(E^{0,1}_\iota)).\ch_1(\omega_Y(S))^{\dim(Y)-1} \leq 0$ \ \ \ \ \ \ 
for all $\iota\in \{1,\ldots, r\}$.\vspace{.1cm}
\item[B.] \hspace*{0.5cm} $\mu(E^{1,0}_1) > \mu(F^{1,0}_2) > \cdots > \mu(F^{1,0}_r) > 0$ \hspace*{\fill} \ \vspace{.1cm}\\
\  \hspace*{\fill} $ 0 > \mu(E^{0,1}_1) >  \mu(F^{0,1}_2) > \cdots > \mu(F^{0,1}_{r})$.\hspace*{0.5cm} \ 
\end{enumerate}
\end{claim}
\begin{proof} By \eqref{eq6.3} $(E_\iota = E^{1,0}_\iota \oplus E^{0,1}_\iota, \theta|_{E_\iota})$ is a Higgs subbundle of $(E,\theta)$. So A) follows from Proposition~\ref{6.3}. Since $(E^{0,1}_1,0)$ is a Higgs subbundle of $(E,\theta)$ and since $(F^{1,0}_r,0)$ is a quotient Higgs bundle, one also obtains $\mu(F^{1,0}_r) > 0 > \mu(E^{0,1}_1)$.

The slope inequalities  
$$
\mu(G_{j_\iota}/G_{j_\iota-1}) > \mu(G_{j_{\iota}+1}/G_{j_{\iota}})\mbox{ \ \ and \ \ } \mu(G'_{j'_\iota}/G'_{j'_\iota-1}) > \mu(G'_{j'_{\iota}+1}/G'_{j'_{\iota}}),
$$ 
together with~\eqref{eq6.4} and~\eqref{eq6.5}, imply the remaining inequalities in B).
\end{proof}
\begin{claim}\label{6.6}\ \hspace*{\fill} $\displaystyle \mu(E^{1,0})-\mu(E^{0,1}) \leq {\rm Max}\{\mu(F^{1,0}_\kappa)-\mu(F^{0,1}_\kappa);\ \kappa= 1,\ldots, r\}$\hspace*{\fill} \ \vspace{.2cm}\\
and the equality is strict except if $r=1$.
\end{claim}
Before proving Claim~\ref{6.6} let us finish the proof of Theorem~\ref{6.4}.
By~\eqref{eq6.3} the Higgs field $\theta$ induces a non-zero map
\begin{equation}\label{eq6.6}
G_{j_{\iota-1} + 1}/G_{j_{\iota-1}} \longrightarrow \big(G'_{j'_{\iota}}/G'_{j'_{\iota}-1}\big) \otimes \Omega^1_Y(\log S).
\end{equation}
The semistability of both sides of~\eqref{eq6.6} implies that
$$
\mu(G_{j_{\iota-1} + 1}/G_{j_{\iota-1}}) \leq \mu(G'_{j'_{\iota}}/G'_{j'_{\iota}-1}) + \mu(\Omega^1_Y(\log S)).
$$
By~\eqref{eq6.4} and~\eqref{eq6.5} one has 
\begin{equation}\label{eq6.7}
\mu(G_{j_{\iota-1}+1}/G_{j_{\iota-1}}) \geq \mu(F^{1,0}_\iota)
\mbox{ \ \ and \ \ }
\mu(F^{0,1}_\iota) \geq \mu(G'_{j'_\iota}/G'_{j'_\iota-1}).
\end{equation}
and altogether
\begin{equation}\label{eq6.8}
 \mu(F^{1,0}_\iota)-\mu(F^{0,1}_\iota) \leq \mu(G_{j_{\iota-1}+1}/G_{j_{\iota-1}}) -
\mu(G'_{j'_\iota}/G'_{j'_\iota-1}) \leq \mu(\Omega^1_Y(\log S)).
\end{equation}
For $j=r$ the first part of Claim~\ref{6.6} implies that $\mu(E^{1,0})-\mu(E^{0,1}) \leq \mu(\Omega^1_Y(\log S))$
as claimed in~\eqref{eq6.2}. This can only be an equality if $r=1$, hence $j_1=\ell$ and $j'_1=\ell'$.

In addition, the equality in~\eqref{eq6.2} can only hold if~\eqref{eq6.8} is an equality. Then 
the two inequalities in~\eqref{eq6.7} have to be equalities as well. By the definition of the Harder-Narasimhan filtration the equalities
$$
\mu(G_{1}) = \mu(E^{1,0})
\mbox{ \ \ and \ \ }
\mu(E^{0,1}) = \mu(G'_{\ell'}/G'_{\ell'-1})
$$
imply that $\ell=\ell'=1$, hence that $E^{1,0}$ and $E^{0,1}$ are both $\mu$-semistable.
\end{proof}
\begin{proof}[Proof of Claim~\ref{6.6}] We will try to argue by induction on the length of the filtration, starting with the trivial case $r=1$. Unfortunately this forces us to replace the
rank of the $F^{0,1}_i$ by some virtual rank. We define:
\begin{enumerate}
\item[(1)] $\gamma_{i}=\ch_1(F_{i}).\ch_1(\omega_Y(S))^{\dim(Y)-1}$.\vspace{.1cm}
\item[(2)] $\mu_{i}^{p,q}= \mu(F_i^{p,q})$ and $\Delta_i=\mu_{i}^{1,0}-\mu_{i}^{0,1}$.\vspace{.1cm}
\item[(3)] $\rho_i^{1,0}=\rk(F_i^{1,0})$ and $\rho_i^{0,1}=\rk(F_i^{0,1})-\ds\frac{\gamma_i}{\mu_i^{0,1}}$.\vspace{.1cm}
\item[(4)] For $0 < \kappa \leq \ell$ 
$$
s_\kappa^{p,q}=\sum_{i=1}^\kappa \rho_i^{p,q}, \ \ \ \ \
\varUpsilon^{1,0}_\kappa=\frac{\sum_{i=1}^\kappa \mu_i^{1,0}\cdot \rho_i^{1,0}}{s_\kappa^{1,0}}, \ \ \ \ \
\varUpsilon^{0,1}_\kappa= \frac{\sum_{i=1}^\kappa \mu_i^{0,1}\cdot \rho_i^{0,1}}{s_\kappa^{0,1}},
$$
and $\delta_\kappa=\varUpsilon_\kappa^{1,0}-\varUpsilon_\kappa^{0,1}$.
\end{enumerate}
Remark that $\varUpsilon^{1,0}_\kappa$ is the slope of the sheaf $E^{1,0}_\kappa$, whereas
$\varUpsilon^{0,1}_\kappa$ is just a {\em virtual slope} without any geometric meaning.

By the choice of $\rho_i^{0,1}$ one finds 
$$
\rho_i^{1,0}\cdot \mu_i^{1,0}+\rho_i^{0,1}\cdot \mu_i^{0,1}= \rk(F_i^{1,0})\cdot \mu_i^{1,0} +
\rk(F_i^{0,1})\cdot \mu_i^{0,1} - \gamma_i =0
$$
and we can state:\vspace{.1cm}
\begin{enumerate}
\item[(5)] $\rho_i^{1,0}\cdot \mu_i^{1,0}= - \rho_i^{0,1}\cdot \mu_i^{0,1}$ and hence $\rho_i^{0,1}>0$.\vspace{.1cm} 
\end{enumerate}
Recall that the condition B) in Claim~\ref{6.4} says that $-\mu_\kappa^{1,0} > -\mu_i^{1,0}$ and
$\mu_i^{0,1} > \mu_\kappa^{0,1}$ for $i<\kappa$. This implies
\begin{multline*}
s_\kappa^{1,0}\cdot \rho_\kappa^{0,1}\cdot \mu_\kappa^{0,1} =\sum_{i=1}^\kappa \rho_i^{1,0}\cdot \rho_\kappa^{0,1} \cdot \mu_\kappa^{0,1}
= \sum_{i=1}^\kappa \rho_i^{1,0}\cdot \rho_\kappa^{1,0} \cdot (-\mu_\kappa^{1,0}) \geq
\sum_{i=1}^\kappa \rho_i^{1,0}\cdot \rho_\kappa^{1,0} \cdot (-\mu_i^{1,0}) =\\
\sum_{i=1}^\kappa \rho_i^{0,1}\cdot \rho_\kappa^{1,0} \cdot \mu_i^{0,1} \geq
\sum_{i=1}^\kappa \rho_i^{0,1}\cdot \rho_\kappa^{1,0} \cdot \mu_\kappa^{0,1}=
s_\kappa^{0,1}\cdot \rho_\kappa^{1,0} \cdot \mu_\kappa^{0,1}.
\end{multline*}
Since $\mu_\kappa^{0,1}$ is negative, one gets\vspace{.1cm}
\begin{enumerate}
 \item[(6)] $s_\kappa^{1,0}\cdot \rho_\kappa^{0,1}\leq  s_\kappa^{0,1}\cdot \rho_\kappa^{1,0}$ or equivalently
$s_{\kappa-1}^{1,0}\cdot \rho_\kappa^{0,1}\leq  s_{\kappa-1}^{0,1}\cdot \rho_\kappa^{1,0}$.\vspace{.1cm} 
\end{enumerate}
The induction step will use the next claim.
\begin{claim}\label{6.7}
For $0 < \kappa \leq \ell$ one has $\delta_\kappa \leq {\rm Max}\{ \delta_{\kappa-1}, \ \Delta_\kappa \},$
with equality if and only if $\delta_{\kappa-1}=\Delta_\kappa$ and $\rho^{1,0}_\kappa\cdot s_{\kappa}^{0,1}=\rho_\kappa^{0,1}\cdot s_{\kappa}^{1,0}$.
\end{claim}
\begin{proof}
We let $A=s_{\kappa-1}^{1,0}\cdot s_{\kappa-1}^{0,1}$, $B=\rho_\kappa^{1,0}\cdot \rho_\kappa^{0,1}$, 
$C=s_{\kappa-1}^{1,0}\cdot \rho_\kappa^{0,1}$ and $D=\rho_\kappa^{1,0}\cdot s_{\kappa-1}^{0,1}$. 
By (6) one has $D-C \geq 0$. Then
\begin{multline*}
s_\kappa^{1,0}\cdot s_\kappa^{0,1}\cdot \delta_\kappa =
\sum_{i=1}^\kappa\big( \mu_i^{1,0}\cdot \rho_i^{1,0}\cdot s_\kappa^{0,1}
- \mu_i^{0,1}\cdot \rho_i^{0,1}\cdot s_\kappa^{1,0} \big)=\\
\mu_\kappa^{1,0}\cdot \rho_\kappa^{1,0}\cdot s_\kappa^{0,1}
- \mu_\kappa^{0,1}\cdot \rho_\kappa^{0,1}\cdot s_\kappa^{1,0} + \sum_{i=1}^{\kappa-1}\big( \mu_i^{1,0}\cdot \rho_i^{1,0}\cdot s_\kappa^{0,1}
- \mu_i^{0,1}\cdot \rho_i^{0,1}\cdot s_\kappa^{1,0} \big)=\\
B\cdot \Delta_\kappa + A\cdot \delta_{\kappa-1} + C\cdot(\varUpsilon_{\kappa-1}^{1,0} - \mu_\kappa^{0,1})
+ D\cdot(\mu_\kappa^{1,0} - \varUpsilon_{\kappa-1}^{0,1})=\\
B\cdot \Delta_\kappa + A\cdot \delta_{\kappa-1} + C\cdot(\delta_{\kappa-1}+\Delta_\kappa)
+ (D-C) \cdot(\mu_\kappa^{1,0} - \varUpsilon_{\kappa-1}^{0,1}).
\end{multline*}
Since $\mu_\kappa^{1,0}  < \mu_i^{1,0}$ for $i < \kappa$ one finds $\mu_\kappa^{1,0} < \varUpsilon_{\kappa-1}^{1,0}$
and 
$$ 
(A+B+C+D) \cdot \delta_\kappa \leq B\cdot \Delta_\kappa + A\cdot \delta_{\kappa-1} + C\cdot\Delta_\kappa
+ D \cdot\delta_{\kappa-1}.
$$
This implies the inequality in Claim~\ref{6.7}. If the equality holds, $\Delta_\kappa=\delta_{\kappa-1}$ and
$$0=D-C=\rho_\kappa^{1,0}\cdot s_{\kappa-1}^{0,1}- s_{\kappa-1}^{1,0}\cdot \rho_\kappa^{0,1}=
\rho_\kappa^{1,0}\cdot s_{\kappa}^{0,1}- s_{\kappa}^{1,0}\cdot \rho_\kappa^{0,1}.
$$
\end{proof}
\begin{claim}\label{6.8}
One has the inequality $\mu(E^{1,0})-\mu(E^{0,1}) \leq \delta_r$ and the equality can only hold for
$\gamma_1=\cdots =\gamma_r=0$.
\end{claim}
\begin{proof}
Since $\mu(E^{1,0})=\varUpsilon_r^{1,0}$ it remains to verify that $\mu(E^{0,1}) \geq \varUpsilon_r^{0,1}$.
As a first step,
\begin{multline}\label{eq6.9}
\Big(\sum_{i=1}^r \rho_i^{0,1}\Big) -\rk(E^{0,1})= \sum_{i=1}^r \big(\rho_i^{0,1}-\rk(F_i^{0,1})\big)=
\sum_{i=1}^r \frac{-\gamma_i}{\mu_i^{0,1}}=\\
\frac{-\gamma_i}{\mu_r^{0,1}} \cdot \Big( \sum_{i=1}^r \gamma_i\big) + \sum_{i=1}^{r-1} \frac{\mu_i^{0,1}-\mu_{i+1}^{0,1}}{\mu_i^{0,1}\cdot \mu_{i+1}^{0,1}}\cdot \Big( \sum_{j=1}^i \gamma_j\big). 
\end{multline}
Since $\sum_{j=1}^i \gamma_j\leq 0$ and equal to zero for $i=r$, and since $\frac{\mu_i^{0,1}-\mu_{i+1}^{0,1}}{\mu_i^{0,1}\cdot \mu_{i+1}^{0,1}}$ is positive, one obtains
$$
\sum_{i=1}^r \rho_i^{0,1} \leq \rk(E^{0,1}).
$$
Then
\begin{multline*}
\mu(E^{0,1}) = \frac{\sum_{i=1}^r \mu_i^{0,1}\cdot \rk(F_i^{0,1})}{\rk(E^{0,1})}=
\frac{\sum_{i=1}^r \mu_i^{0,1}\cdot \rho_i^{0,1}}{\rk(E^{0,1})}+
\frac{\sum_{i=1}^r \gamma_i^{0,1}}{\rk(E^{0,1})}=\\
\frac{\sum_{i=1}^r \mu_i^{0,1}\cdot \rho_i^{0,1}}{\rk(E^{0,1})}
\geq \frac{\sum_{i=1}^r \mu_i^{0,1}\cdot \rho_i^{0,1}}{\sum_{i=1}^r \rho_i^{0,1}}=\varUpsilon_r^{0,1},  
\end{multline*}
as claimed. The equality implies that the expression in~\eqref{eq6.9} is zero, which is only possible if $\gamma_1=\cdots =\gamma_r=0$.
\end{proof}
Using the Claims~\ref{6.7} and~\ref{6.8} one finds that
\begin{multline*}
\mu(E^{1,0})-\mu(E^{0,1}) \leq  \delta_r \leq {\rm Max}\{ \delta_{r-1},\Delta_r\} \leq
{\rm Max}\{ \delta_{r-2},\Delta_{r-1}, \Delta_r\}\leq\cdots \\
\cdots \leq {\rm Max}\{ \Delta_1, \ldots ,\Delta_{r-1}, \Delta_r\}.
\end{multline*}
The equality implies that for all $\kappa$ the inequalities in Claims~\ref{6.7} and~\ref{6.8} are equalities.
The second one implies that for all $\kappa$ one has $\gamma_\kappa=0$, hence $\rho_\kappa^{0,1}=\rk(F_\kappa^{0,1})$, and the first one that 
$$
0=\rho^{1,0}_\kappa\cdot s_{\kappa}^{0,1}-\rho_\kappa^{0,1}\cdot s_{\kappa}^{1,0}=
\rk(F^{1,0}_\kappa)\cdot s_{\kappa}^{0,1}-\rk(F^{0,1}_\kappa)\cdot s_{\kappa}^{1,0}.
$$
\end{proof}
As for variation of Hodge structures over curves, the Arakelov inequality \eqref{6.2} is a direct consequence of the polystability of the Higgs bundle $(E,\theta)$. The Arakelov equality $\mu(\V) = \mu(\Omega^1_Y(\log S))$ allows to deduce the semistability of the sheaves $E^{1,0}$ and $E^{0,1}$. 
However, we do not know whether one gets the stability, as it has been the case over curves (see \ref{3.4}). 
Although we were unable to construct an example, we do not expect this.  

So it seems reasonable to ask, which additional conditions imply the stability of the sheaves
$E^{1,0}$ and $E^{0,1}$.  

\section{Geodecity of higher dimensional subvarieties in $\sA_g$}\label{s.7}
Let us recall the geometric interpretation of the Arakelov equality, shown in
\cite{vz07} and \cite{mvz07}. 
\begin{assumptions}\label{7.1}
We keep the assumptions and notations from Section~\ref{s.6}. Hence $Y$ is a projective non-singular
manifold, and $Y_0\subset Y$ is open with $S=Y\setminus Y_0$ a normal crossing divisor. We assume the positivity condition ($\star$) and we consider an irreducible polarized $\C$-variation of Hodge structures $\V$ of weight one with unipotent monodromies around the components of $S$. As usual its Higgs bundle will be denoted by $(E,\theta)$.
\end{assumptions}
The first part of Yau's Uniformization Theorem (\cite{Ya93}, discussed in \cite[Theorem 1.4]{vz07}) was already used in the last section. It says that the Assumption ($\star$) forces the sheaf $\Omega_Y^1(\log S)$ to be $\mu$-polystable. The second part gives a geometric interpretation of stability properties of the direct factors. Writing
\begin{equation}\label{eq7.1}
\Omega_Y^1(\log S)=\Omega_1\oplus \cdots \oplus \Omega_s
\end{equation}
for its decomposition as direct sum of $\mu$-stable sheaves and $n_i=\rk(\Omega_i)$, 
we say that {\em $\Omega_i$ is of type A}, if it is invertible, 
and {\em of type B}, 
if $n_i>1$ and if for all $\ell>0$ the sheaf $S^\ell(\Omega_i)$ is $\mu$-stable. 
In the remaining cases, i.e.\ if 
for some $\ell>1$ the sheaf $S^\ell(\Omega_i)$ is $\mu$-unstable, 
we say that $\Omega_i$ is of type C.

Let $\pi:\tilde{Y}_0 \to Y_0$ denote the universal covering with  covering group $\Gamma$. The decomposition (\ref{eq7.1}) of $\Omega^1_Y(\log S)$ gives rise to a product structure 
\begin{equation*}
\tilde{Y}_0=M_1\times \cdots \times M_s,
\end{equation*}
where $n_i=\dim(M_i)$. The second part of Yau's Uniformization Theorem gives a criterion for each $M_i$ to be a bounded symmetric domain. This is automatically the case if $\Omega_i$ is of type A or C.
If $\Omega_i$ is of type B, then $M_i$ is a $n_i$-dimensional complex ball if and only if
\begin{equation}\label{eq7.2}
\big[2\cdot (n_i+1)\cdot \ch_2(\Omega_i)-n_i\cdot \ch_1(\Omega_i)^2\big].\ch(\omega_Y(S))^{\dim(Y)-2}=0.
\end{equation}
\begin{definition}\label{7.2}
The variation of Hodge structures $\V$ is called pure (of type $i$) if the Higgs field
factors like
$$
E^{1,0} \longrightarrow E^{0,1}\otimes \Omega_i \subset E^{0,1}\otimes \Omega_Y^1(\log S)
$$
(for some $i=i(\V)$).
\end{definition}
If one knows that $\tilde{Y}_0$ is a bounded symmetric domain, hence if 
\eqref{eq7.2} holds for all direct factors of type B, one obtains the purity of $\V$ as a consequence of
the Margulis Superrigidity Theorem:
\begin{theorem}\label{7.3}
Suppose in~\ref{7.1} that $\tilde{Y}_0$ is a bounded symmetric domain. Then $\V$ is pure.  
\end{theorem}
\begin{proof}[Sketch of the proof]
Assume first that $Y_0=U_1\times U_2$. By \cite[Proposition 3.3]{vz05} an irreducible 
local system on $\V$ is of the form ${\rm pr}_1^*\V_1\otimes {\rm pr}_2^*\V_2$,
for irreducible local systems $\V_i$ on $U_i$ with Higgs bundles
$(E_i,\theta_i)$. Since $\V$ is a variation of Hodge structures of weight $1$,
one of those, say $\V_2$ has to have weight zero, hence it must be unitary. 

Then the Higgs field on $Y_0$ factors through $E^{0,1}\otimes \Omega^1_{U_1}$.
By induction on the dimension we may assume that $\V_1$ is pure of type $\iota$
for some $\iota$ with $M_\iota$ a factor of $\tilde{U}_1$. Hence the same holds true for
$\V$.

So we may assume that all finite \'etale coverings of $Y_0$ are indecomposable.
By \cite{Zi} \S~2.2, replacing $\Gamma$ by a subgroup of finite index, hence replacing $Y_0$ by a finite unramified cover, there is a partition of $\{1,\ldots,s\}$ into subsets $I_k$ such that $\Gamma = \prod_k \Gamma_k$ and $\Gamma_k$ is an irreducible lattice in $\prod_{i \in I_k} G_i$. Here irreducible
means that for any normal subgroup $N \subset \prod_{i \in I_k} G_i$ the image of
$\Gamma_k$ in $\prod_{i \in I_k} G_i/N$ is dense. 
Since the finite \'etale coverings of $Y_0$ are indecomposable, $\Gamma$ is irreducible,
so $I_1 =\{1,\ldots,s\}$. 
\par
If $s=1$ or if $\V$ is unitary, the statement of the proposition is trivial. 
Otherwise, $G:=\prod_{i=1}^s G_i$ is of real rank $\geq 2$ and the conditions
of Margulis' superrigidity theorem (e.g.\ \cite[Theorem 5.1.2 ii)]{Zi})
are met.  As consequence, the homomorphism $\Gamma \to \Sp(V,Q)$,
where $V$ is a fibre of $\V$ and where $Q$ is the symplectic form on $V$,
factors through a representation $\rho: G \to \Sp(V,Q)$. Since the 
$G_i$ are simple, we can repeat the argument from \cite[Proposition 3.3]{vz05}, used above
in the product case: 
$\rho$ is a tensor product of representations, all of which but one have weight $0$. 
\end{proof}
The next theorem replaces the condition that $\tilde{Y}_0$ is a bounded symmetric domain by the Arakelov equality.
\begin{theorem}\label{7.4}
Suppose in~\ref{7.1} that $\V$ satisfies the Arakelov equality 
$$
\mu(\V) = \mu(\Omega^1_Y(\log S)).
$$
Then $\V$ is pure.
\end{theorem}
The two Theorems~\ref{6.4} and~\ref{7.4} imply that the Higgs field of $\V$ is given by a morphism
$$
E^{1,0} \longrightarrow E^{0,1}\otimes \Omega_i
$$
between $\mu$-semistable sheaves of the same slope. If $\Omega_i$ is of type A or C this implies geodecity (for the Hodge or Bergman metric) in period domains of variation of Hodge structures of weight one.
\begin{theorem}\label{7.5} Suppose in Theorem~\ref{7.4} that for $i=i(\V)$ the sheaf $\Omega_i$ is of type A or C.
Let $M'$ denote the period domain for $\V$. Then the period map factors as the projection $\tilde{Y}_0 \to M_i$ and a totally geodesic embedding $M_i \to M'$.
\end{theorem} 
If $\Omega_i$ is of type B we need some additional numerical invariants in order to deduce a similar property.

Let $(F,\tau)$ be any Higgs bundle, not necessarily of degree zero. For $\ell=\rk(F^{1,0})$ consider the Higgs bundle
\begin{gather}\notag
\bigwedge^\ell (F,\tau)=\big(\bigoplus_{i=0}^\ell F^{\ell-i,i}, \ \bigoplus_{i=0}^{\ell-1} \tau_{\ell-i,i} \big)
\mbox{ \ \ with}\\ \label{eq7.3}
F^{\ell-m,m}=\bigwedge^{\ell-m}(F^{1,0}) \otimes \bigwedge^m (F^{0,1})\mbox{ \ \ and with}\\ \notag
\tau_{\ell-m,m}: \bigwedge^{\ell-m}(F^{1,0}) \otimes \bigwedge^m (F^{0,1}) \longrightarrow
\bigwedge^{\ell-m-1}(F^{1,0}) \otimes \bigwedge^{m+1} (F^{0,1}) \otimes \Omega_Y^1(\log S)
\end{gather}
induced by $\tau$. Then  $F^{\ell,0}=\det(F^{1,0})$ and $\langle\det(F^{1,0})\rangle$ denotes the Higgs subbundle of $\bigwedge^\ell(F,\tau)$ generated by $\det(F^{1,0})$. Writing 
$$
\tau^{(m)}=\tau_{\ell-m+1,m-1}\circ \cdots \circ \tau_{\ell,0},
$$
we define as a measure for the complexity of the Higgs field
\begin{multline*}
\varsigma((F,\tau)): = {\rm Max}\{\ m\in \N ; \ \tau^{(m)}(\det(F^{1,0}))\neq 0\}=\\
{\rm Max}\{\ m\in \N ; \ \langle\det(F^{1,0})\rangle^{\ell-m,m} \neq 0\}.
\end{multline*}
For the Higgs bundle $(E,\theta)$ of $\V$, we write
$\varsigma(\V)=\varsigma((E,\theta))$.
\begin{lemma}\label{7.6}
Suppose in~\ref{7.1} that $\V$ satisfies the Arakelov equality and, using the notation from 
Theorem~\ref{7.4}, that for $i=i(\V)$ the sheaf $\Omega_i$ is of type B (or of type A). Then 
\begin{equation}\label{eq7.4}
\varsigma(\V) \geq \frac{\rk(E^{1,0})\cdot\rk(E^{0,1})\cdot(n_i+1)}{\rk(E) \cdot n_i}.
\end{equation}
Moreover \eqref{eq7.4} is an equality if and only if the kernel of the morphism 
$$
\sH om(E^{0,1},E^{1,0}) \to \Omega^1_Y(\log S),
$$ 
induced by $\theta$, is a direct factor of $\sH om(E^{0,1},E^{1,0})$. 
\end{lemma}
Here again one uses Simpson's polystability, applied to the variation of Hodge structures \ \  $\displaystyle \bigwedge^\ell \V$ \ \ with Higgs  bundle \ \ $\displaystyle \bigwedge^\ell (E,\theta)$.

\begin{theorem}\label{7.7} Suppose in Theorem~\ref{7.4} that for $i=i(\V)$ the sheaf $\Omega_i$ is of type A or B.
Assume that one has the length equality
\begin{equation}\label{eq7.5}
\varsigma(\V) = \frac{\rk(E^{1,0})\cdot\rk(E^{0,1})\cdot(n_i+1)}{\rk(E) \cdot n_i}.
\end{equation}
Then
\begin{enumerate}
\item[a.] $M_i$ is the complex ball $\SU(1,n_i)/K$, and $\V$ is the tensor product of a unitary representation with a wedge product of the standard representation of $\SU(1,n_i)$.
\item[c.] Let $M'$ denote the period domain for $\V$. Then the period map factors as the projection $\tilde{Y}_0 \to M_i$ and a totally geodesic embedding $M_i \to M'$.
\end{enumerate}
\end{theorem} 
In Theorem~\ref{7.7}, a) the Higgs field of the standard representation of $\SU(1,n_i)$ (or of its dual)
is given by
$$
E^{1,0}=\omega_i^{-\frac{1}{n_i+1}}\otimes \Omega_i, \ \ \ 
E^{0,1}=\omega_i^{-\frac{1}{n_i+1}} \mbox{ \ \ and \ \ }
\theta={\rm id}: \omega_i^{-\frac{1}{n_i+1}}\otimes \Omega_i \longrightarrow
\omega_i^{-\frac{1}{n_i+1}}\otimes \Omega_i,
$$ 
where $\omega_i^{-\frac{1}{n_i+1}}$ stands for an invertible sheaf, whose 
$(n_i+1)$-st power is
$\det(\Omega_i)$. 

\begin{remark}\label{7.8}
We do not know, whether the Arakelov equality implies the condition~\eqref{eq7.5}. 
In \cite{mvz07} this implication has been verified for $\rk(\V)\leq 7$. Nevertheless, the necessity of the Yau-equality in the characterization of complex ball quotients indicates that besides of the Arakelov equality one needs a second condition, presumably one using second Chern classes. 
\end{remark}
Although the second Chern class does not occur in Theorem \ref{7.7} it seem to be hidden in the condition on the length of the Higgs field stated there. As an illustration of the latter, let us consider a second 
numerical invariant, the discriminant. Recall that for a torsion free coherent sheaf $\sF$ on $Y$
$$
\delta(\sF)=\big[2\cdot \rk(\sF)\cdot\ch_2(\sF)-(\rk(\sF)-1)\cdot \ch_1
(\sF)^2\big].\ch_1(\omega_Y(S))^{\dim(Y)-2}.
$$
For the Higgs bundle $(E,\theta)$ of $\V$ we define
$\delta(\V)= {\rm Min}\{\delta(E^{1,0}), \delta(E^{0,1})\}$.

The Bogomolov inequality for semi-stable locally free sheaves allows to state as a corollary of Theorem~\ref{6.4}:
\begin{corollary}\label{7.9}
Keeping the assumptions and notations from Theorem~\ref{6.4} the Arakelov equality
$\mu(\V)= \mu (\Omega^1_Y(\log S))$ implies that $\delta(\V)\geq 0$.
\end{corollary}
\begin{theorem}\label{7.10} Suppose in Theorem~\ref{7.4} that $\omega_Y(S)$ is ample, that for $i=i(\V)$ the sheaf $\Omega_i$ is of type A or B and that $\delta(\V)=0$. Then 
\begin{enumerate}
\item[a.] $M_i$ is the complex ball $\SU(1,n_i)/K$, and $\V$ is the tensor product of a unitary representation with the standard representation of $\SU(1,n_i)$.
\item[b.] Let $M'$ denote the period domain for $\V$. Then the period map factors as the projection $\tilde{Y}_0 \to M_i$ and a totally geodesic embedding $M_i \to M'$.\vspace{.1cm}
\item[c.] \ \hspace*{\fill} $\ds
\varsigma(\V) = \frac{\rk(E^{1,0})\cdot\rk(E^{0,1})\cdot(n_i+1)}{\rk(E) \cdot n_i}.
$ \hspace*{\fill} \ 
\end{enumerate}
\end{theorem} 
Note that in a) we have to exclude the wedge products of the standard representations. For those
$\delta(\V)$ is larger than $0$.

In \cite{mvz07} we are mainly interested in subvarieties of $\sA_g$. If one assumes the conditions 
in~\ref{7.1} to hold for all non-unitary local $\C$-subvariations of Hodge structures of the induced family
then one can deduce the following numerically characterization of Kuga fibre spaces: 
\begin{theorem}\label{7.11}
Let $f:A \to Y_0$ be a family of polarized abelian varieties such that
$R^1f_{0*}\C_A$ has unipotent local monodromies at infinity, and such that the induced 
morphism $Y_0\to \sA_g$ is generically finite. Assume that $Y_0$ has a projective compactification $Y$ satisfying the Assumption ($\star$).

Then the following two conditions are equivalent:
\begin{enumerate}
\item[I.] There exists an \'etale covering $Y_0'\to Y_0$ such that the pullback family
$f':A'=A\times_{Y_0}Y_0' \to Y_0'$ is a Kuga fibre space.
\item[II.] For each irreducible subvariation of Hodge structures $\V$ of $R^1f_{*0}\C_A$
with Higgs bundle $(E,\theta)$ one has:
\begin{enumerate}
\item[1.] Either $\V$ is unitary or the Arakelov equality $\mu(\V)=\mu(\Omega_Y^1(\log S))$ holds.
\item[2.] If for a $\mu$-stable direct factor $\Omega_j$ of $\Omega^1_Y(\log S)$ 
of type B the composition
$$
\theta_j:E^{1,0} \stackrel{ \theta}{\longrightarrow} E^{0,1}\otimes \Omega_Y^1(\log S)  \stackrel{ {\rm pr}}{\longrightarrow} E^{0,1}\otimes \Omega_j
$$
is non-zero, then 
$$
\varsigma((E,\theta_j)) = \frac{\rk(E^{1,0})\cdot\rk(E^{0,1})\cdot(n_j+1)}{\rk(E) \cdot n_j}.
$$
\end{enumerate}
\end{enumerate}
\end{theorem}
\begin{remarks} \ \\[.1cm]
(1) Theorem \ref{7.11} partly answers the question on the $\mu$-stability of the Hodge bundles $E^{1,0}$
and $E^{0,1}$, at least for subvariations of Hodge structures in $R^1f_{0*}\C_{X_0}$ for a family $f_0:X_0\to Y_0$ of abelian varieties. In fact, choosing in part I) a Mumford compactification $Y'$ of $Y_0'$ one can show that
the Hodge sheaves $E'^{1,0}$ and $E'^{0,1}$ of the pullback $\V'$ of the irreducible subvariation
of Hodge structures $\V$ are $\mu$-stable. So up to replacing $Y_0$ by an \'etale cover and $Y$ by a suitable
compactification, the $\mu$-stability of the Hodge sheaves follows from the Arakelov equality if
$\V$ is of type A or C, whereas for type B we need an additional numerical condition.\\[.1cm]
(2) If one knows the $\mu$-stability of  $E'^{1,0}$ and $E'^{0,1}$ on some compactification of an \'etale covering $Y_0'$ of $Y_0$, and if $\omega_Y(S)$ is ample, then $\sH om(E^{0,1},E^{1,0})$ is $\mu$-polystable
and the Arakelov equality implies that the morphism 
$$
\sH om(E^{0,1},E^{1,0}) \to \Omega^1_Y(\log S),
$$ 
is surjective and splits.
So by Lemma \ref{7.6} the numerical condition, saying that \eqref{eq7.4} is an equality, holds and by Theorem
\ref{7.7} $M_i$ must be a complex ball. As remarked in \ref{7.8} we think it is unlikely to have a characterization of a complex ball, which is only using first Chern classes. \\[.1cm]
(3) The condition ``$\omega_Y(S)$ ample'' appears in (2) since one uses that the tensor product of polystable sheaves is polystable. The same is used in the proof of Theorem \ref{7.10}. There however the ampleness is needed for a second reason. One uses the characterization of unitary bundles as those polystable bundles with vanishing first and second Chern class. S.T. Yau conjectures that for both statements ``$\omega_Y(S)$ nef and big'' is sufficient. He and Sun promised to work out a proof of those results.  
\end{remarks} 
\section{Open ends}\label{s.8}
{\bf I.} As mentioned already, under the assumptions made in \ref{7.1} for variations of Hodge structures 
$\V$ of weight one and of small rank, the Arakelov equality implies that the length inequality \ref{eq7.4} is
an equality. Let us write in \ref{7.1} $q=\rk(E^{1,0})$, $p=\rk(E^{0,1})$ and assume that $q\leq p$.
Since $\deg(E^{1,0})+\deg(E^{0,1})=0$ one can rewrite the Arakelov inequality \eqref{6.2} as
\begin{equation}\label{eq8.1}
\ch_1(E^{1,0}).\ch_1(\omega_Y(S))^{\dim(Y)-1} \leq \frac{p\cdot q}{(p+q)\cdot \dim(Y)}\cdot \ch_1(\omega_Y(S))^{\dim(Y)}. 
\end{equation}
\begin{lemma}\label{8.1} Assume that $\V$ satisfies the Arakelov equality and that for $i=i(\V)$ the sheaf $\Omega_i$ is of type A or B. Then 
$$
n_i \cdot q \geq p\geq q, \mbox{ \ \ for \ \ } n_i=\rk(\Omega_i), 
$$ 
and if $p= n_i \cdot q$ the numerical condition \eqref{eq7.5} in Theorem \ref{7.7} holds.
In particular $M_i$ is a complex ball in this case.
\end{lemma}
\begin{proof}
This follows from the definition of $\varsigma((E,\theta))$ and Lemma \ref{7.6}, implying that
\begin{equation}\label{eq8.2}
q \geq \varsigma((E,\theta)) \geq \frac{p\cdot q \cdot (n_i+1)}{(p+q)\cdot n_i}.
\end{equation}  
\end{proof}
Assume the Arakelov equality. If $q=1$ $E^{1,0}$ is invertible, $E^{1,0}\otimes T_Y(-\log S)$ and $E^{0,1}$ have to be $\mu$-equivalent. So $p=m$ and \eqref{eq8.2} must be an equality, as required in \eqref{eq7.5}. Hence $M_i$ is a complex ball of dimension $n_i$ (see \cite[Example 8.5]{mvz07}). 

If $q=2$, assuming that $\omega_Y(S)$ is ample, one can apply \cite[Lemma 8.6 and Example 8.7]{mvz07}, and again one finds that
the Arakelov equality implies the length equality \eqref{eq7.5}.
\begin{corollary}\label{8.2}
Assume in \ref{7.1} that $\V$ satisfies the Arakelov equality, and that for $i=i(\V)$ the sheaf $\Omega_i$ is of type A or B. Assume that 
$$
{\rm Min}\{\rk(E^{1,0}), \rk(E^{0,1})\} \leq 2,  
$$ 
Then $M_i$ is the complex ball $\SU(1,n_i)/K$, and $\V$ is the tensor product of a unitary representation with a wedge product of the standard representation of $\SU(1,n_i)$.
\end{corollary} 
{\bf II.} In \cite{km08a} Koziarz and Maubon define a Toledo invariant for representations $\rho$ of the fundamental group of a projective variety $Y$ with values in certain groups, in particular in ${\rm SU}(q,p)$. They assume that 
$X$ is of general type, and they use the existence of the canonical model $X_{\rm can}$ of $X$, shown in \cite{BCHM}.
Let us assume here for simplicity, that $X$ is the canonical model, hence that $\omega_Y$ is ample.

As in Section \ref{s.5} the Higgs bundle corresponding to $\rho$ is of the form $\sV\oplus \sW$ and the Higgs field has two components
$$
\beta: \sW \longrightarrow \sV \otimes \Omega_Y^1 \mbox{ \ \ and \ \ }
\gamma: \sV \longrightarrow \sW \otimes \Omega_Y^1.
$$
In \cite[Section 4.1]{km08b} the Toledo invariant is identified with $\deg(\sV)=-\deg(\sW)$, and 
for $1\leq q \leq 2 \leq p$ the generalized Milnor-Wood inequalities in \cite[Theorem 3.3]{km08a} and \cite[Proposition 4.3]{km08b} say that
\begin{equation}\label{eq8.3}
|\ch_1(\sV).\ch_1(\omega_Y)^{\dim(Y)-1}| \leq \frac{q}{\dim(Y)+1}\cdot \ch_1(\omega_Y)^{\dim(Y)}. 
\end{equation}
For $q=2$ one finds in \cite[Proposition 1.2]{km08a} a second inequality, saying
\begin{equation}\label{eq8.4}
|\ch_1(\sV).\ch_1(\omega_Y)^{\dim(Y)-1}| \leq \frac{2p}{(p+2)\cdot\dim(Y)}\cdot \ch_1(\omega_Y)^{\dim(Y)}. 
\end{equation}
In \cite[Theorem 4.1]{km08b} the authors also study the case that the Milnor-Wood inequality \eqref{eq8.3} is an equality. They show that this can only happen if $p\geq m\cdot q$, and that the universal covering $\tilde{Y}$ is a complex ball. 

As in Example \ref{5.2} one can apply \eqref{eq8.3} and \eqref{eq8.4} to a polarized variation of Hodge structures of weight one over $Y=Y_0$. So we will assume that $q=\rk(E^{1,0})$ is smaller than or equal to
$p=\rk(E^{0,1})$ and we will write $n=\dim(Y)$. Here the second inequality \eqref{eq8.4} coincides with the Arakelov inequality \eqref{eq8.1}. As pointed out in \cite[Section 3.3.1]{km08a}, for variations of Hodge structures of weight one \eqref{eq8.3} also holds for $q > 2$. 
\begin{proposition}\label{8.3}
In \ref{7.1} one has the Milnor-Wood type inequality
\begin{equation}\label{eq8.5}
(1+n)\cdot \mu(E^{1,0}) \leq n\cdot \mu(\Omega^1_Y(\log S)).
\end{equation}
The equality implies that $p=q\cdot n$ and hence that \eqref{eq8.5}
coincides with the Arakelov (in)equality.

If $S^\nu(\Omega^1_Y(\log S))$ is stable for all $\nu >0$, and if \eqref{eq8.5} is an equality, then the universal covering $M$ of $U$ is the complex ball $\SU(1,n)/K$, and $\V$ is the tensor product of a unitary representation with the standard representation of $\SU(1,n)$. 
\end{proposition}
\begin{proof}
Let us repeat the argument used in \cite{km08a} in the special case of a variation of Hodge structures of weight one, allowing logarithmic poles of the Higgs bundles along the normal crossing divisor $S$.
As in the proof of Theorem \ref{6.4} one starts with the maximal destabilizing $\mu$-semistable subsheaf 
$G$ of $E^{1,0}$. Let $G'$ be the image of $G\otimes T_Y(-\log S)$ in $E^{0,1}$. Then
the $\mu$-semistability of $G\otimes T_Y(-\log S)$ and the choice of $G$ imply
\begin{gather}\label{eq8.6}
\mu(G') \geq \mu(G) - \mu(\Omega^1_Y(\log S)), \ \ \ \  \mu(G) \geq \mu(E^{1,0}),\\ \label{eq8.7}
\mbox{and \ \ } \rk(G')\leq \rk(G)\cdot n.
\end{gather}
Since $(G\oplus G', \theta|_{G\oplus G'})$ is a Higgs subbundle of $(E,\theta)$ one finds
\begin{multline*}
0 \geq \deg(G)+\deg(G')= \rk(G)\cdot \mu(G) + \rk(G') \cdot \mu(G') \geq\\ 
(\rk(G) + \rk(G'))\cdot \mu(G) - \rk(G')\cdot \mu(\Omega^1_Y(\log S)),\\
\mbox{hence \ \ } \mu(\Omega^1_Y(\log S)) \geq \big(1+\frac{\rk(G)}{\rk(G')}\big)\cdot \mu(G)
\geq \big(1+\frac{1}{n}\big)\cdot \mu(E^{1,0}),\hspace*{1cm}
\end{multline*}
as claimed. If $n\cdot \mu(\Omega^1_Y(\log S)) = (1+n)\cdot \mu(E^{1,0})$ one finds that
all the inequalities in \eqref{eq8.6} and \eqref{eq8.7} are equalities. The first one and the irreducibility of $\V$ imply that $G=E^{1,0}$ and that $G'=E^{0,1}$, whereas the last one shows that
$p=n\cdot q$ for $q=\rk(E^{1,0})$ and $p=\rk(E^{0,1})$. Then
$$
\frac{p\cdot q}{(p+q)\cdot n} = \frac{q}{n+1}
$$
and the equality is the same as the Arakelov equality. 

Finally Lemma \ref{8.1} allows to apply Theorem \ref{7.7}, in case that $\Omega^1_Y(\log S)$ is $\mu$-stable and of type A or B. 
\end{proof}
The situation considered in \cite{km08a} and \cite{km08b} is by far more general than the one studied in Proposition \ref{8.3}. Nevertheless the comparism of the inequalities \eqref{eq8.3} and \eqref{eq8.4} seems to indicate that an optimal Milnor-Wood inequality for for representations in  ${\rm SU}(q,p)$ with $q,p >2$ should have a slightly different shape. As said in Remark \ref{7.8}, it is likely that an interpretation of the equality will depend on a second numerical condition.\vspace{.2cm}

{\bf III.} The proof of the Arakelov inequality \eqref{3.4} for $k>1$ and the interpretation of equality break down if the rank of $\Omega^1_Y(\log S)$ is larger than one. In the proof of Theorem \ref{6.4} we used in an essential way that the weight of the variation of Hodge structures is one. For the Milnor-Wood inequality for a representation of the fundamental group of a higher dimensional manifold of general type with values in ${\rm SU}(p,q)$
one has to assume that ${\rm Min}\{p,q\}\leq 2$, which excludes any try to handle variations of Hodge structures
of weight $k>1$ in a similar way we did in Example \ref{5.3}.
So none of the known methods give any hope for a generalizations of the Arakelov inequality to variations of Hodge structures of weight $k>1$ over a higher dimensional base. We do not even have a candidate for an Arakelov inequality. 

On the other hand, in the two known cases the inequalities are derived from the polystability of the Higgs bundles and the Arakelov equalities are equivalent to the Arakelov condition, defined in \ref{2.2}, iii).
So for weight $k>1$ over a higher dimensional base one should try to work directly in this set-up. 

Even for $k>1$ and $\dim(Y)=1$, as discussed in Section \ref{s.3}, we do not really understand the geometric implications of the Arakelov equality \eqref{eq3.5}, even less the possible implications of the Arakelov condition over a higher dimensional base. 
Roughly speaking, the Addendum \ref{3.6} says that the irreducible subvariations of Hodge structures of weight $k$ 
over a curve, which satisfy the Arakelov equality, look like subvariations of the variation of Hodge structures of weight $k$ for a family of $k$-dimensional abelian varieties. However we do not see a geometric construction relating the two sides.\vspace{.2cm}
 
{\bf IV.} Can one extend the results of \cite{mv08}, recalled in Section \ref{s.4}, to higher dimensional bases?
For example, assume that $\overline{\sA}_g$ is a Mumford compactification of a fine moduli scheme $\sA_g$ with a suitable level structure and that $\varphi: Y\to \overline{\sA}_g$ is an embedding. Writing $S_{\overline{\sA}_g}$ for the boundary, assume that $(Y,S=\varphi^{-1}(S_{\overline{\sA}_g}))$ satisfies the condition ($\star$) in Assumption \ref{6.1}. So one would like to characterize the splitting of the tangent map 
$$
T_Y(-\log S) \longrightarrow \varphi^*T_{\overline{\sA}_g}(-\log  S_{\overline{\sA}_g}) 
$$
in terms of the induced variation of Hodge structures, or in terms of geodecity of
$Y$ in $\overline{\sA}_g$.


\end{document}